\documentclass[10pt]{article}

\usepackage{graphicx}
\usepackage{amsmath,amsfonts,amssymb,latexsym,graphicx}
\usepackage{hyperref}
\usepackage{here}

\usepackage{color}
\newtheorem{theo}{Theorem}
\newtheorem{lem}[theo]{Lemma}
\newtheorem{prop}[theo]{Proposition}

\newcommand{\sign}{\mbox{sign}}

\newcommand\Z{\mathbb{Z}}
\newcommand\R{\mathbb{R}}
\newcommand\N{\mathbb{N}}

\newcommand\E{\mathbb{E}}

\title{Identification of the  Multivariate  Fractional Brownian Motion}

\author{Pierre-Olivier Amblard${}^{1,2}$ and  Jean-Fran\c{c}ois Coeurjolly${}^{1,3}$ \\
 ${}^1$ GIPSAlab/CNRS, France\\
${}^2$ Dept Math\&Stat. and Center for Neural Engineering, \\ The University of Melbourne, Parkville, VIC 3010, Australia \\
${}^3$ LJK, UMR 5226, Grenoble University}

\begin{document}
\maketitle

\def\eqloi{\stackrel{\rm fidi}{=}} 

\begin{abstract}
This paper deals with the identification of  the multivariate
fractional Brownian motion, a recently developed extension of the fractional Brownian motion to the multivariate case. This process is a $p$-multivariate self-similar Gaussian process parameterized by $p$ different Hurst exponents $H_i$, $p$ scaling coefficients $\sigma_i$ (of each component) and also by $p(p-1)$ coefficients $\rho_{ij},\eta_{ij}$ (for $i,j=1,\ldots,p$ with $j>i$) allowing two components to be more or less strongly correlated and allowing the process to be time reversible or not. We investigate the use of discrete filtering techniques to estimate jointly or separately the different parameters and prove the efficiency of the methodology with a simulation study and the derivation of asymptotic results.\\

  
\end{abstract}

\section{Introduction, main results}

The last decade has seen a dramatic effort of research to understand real networks, or complex networks, of any kind \cite{Newm10,Spor10,Strog03}. Indeed, many systems whether natural or man-made constitute networks of interacting systems. These networks are usually considered as complex systems, in the sense that a global behavior emerges from the interaction and cannot be predicted from the sole observation of the individuals.  In general,  the complexity of the system gives to measurements taken at individuals difficult properties such has nonstationarity, fractality, long-range dependence, \ldots This for example occurs in functional magnetic resonance imaging (fMRI), where data collected from different parts of the brain are of course correlated between each other, but also present long-range dependence  \cite{AchaSWSB06,AchaBMB08}. In internet tomography, it is now well recognised that time series corresponding to IP packets or bytes are correlated and long-range dependent \cite{AbryBFRV02}. But mutlivariate time series depicting long-range dependence have also been encountered in fields as different as physics or economics \cite{GilAlana03,AriCar09}.

When measurements are collected simultaneously at several nodes of the networks, the global data set has to be modeled as a multivariate time series. Conversely,  given a multivariate signal, a goal may be to solve an inverse problem: identification of the network underlying the multivariate measurement (each component is associated to a node of the network; a  link between two nodes assesses for dependence between the components.) This problem is a problem of graphical modeling \cite{Whit89,Laur96}. To model long-range multivariate processes, we studied in \cite{CoeuAA10} the extension to the multivariate case of the fractional Brownian motion (and its increments). The mfBm is a Gaussian multivariate signal, whose components are correlated scalar fBm with {\it a priori} different Hurst exponents. This model is interesting for modeling fMRI data. In this paper, we work for the converse problem, developing a methodology to identify the mfBm.

The multivariate fractional Brownian motion is characterized by the Hurst exponents of its components, by  its covariance matrix at time 1, and also by an antisymmetric matrix $\eta_{ij}$ which controls the time assymmetry of the multivariate process. We provide here a framework to estimate all these characteristics from the observation of one sample path of the mfBm. This multivariate process is a nonstationary process with stationary increments. Thus in order to perform time average we work on the increments directly. However, as we recall in  section \ref{mfbmfact:ssec}, the components of the increments process may be long-range dependent individually, and may also present what we call long-range interdependence, meaning that their cross-correlation function may be not summable. This leads to considerable difficulties in the inference methods, especially implying very poor convergence rates (see \cite{Chung02} for example). To circumvent the problem in the scalar case, it is well-known that derivatives smoother than increments have to be considered. The most popular smooth derivative is provided by the wavelet transform when the wavelet is chosen to be orthogonal to polynomial of sufficient high degrees  (see \cite{Flan89,Flan92,TewfK92} for early references). Here, we use a slightly different approach using discrete, compactly supported filters, that need to be orthogonal to some polynomial, but are not necessarily linked to wavelet theory (in that they do not necessarily are the base for a multiresolution analysis).

The filtering is performed for  dilated version of the filter with factor $m$. Each component of the multivariate signal is so filtered. 
 We show that the cross-covariance between the components of the filtered version is a power law of $m$. 
 This generalizes the well-known power law behavior as a function of scale of the variance of the wavelet coefficients in the scalar fractional Brownian motion case.  
 Thus we perform a linear regression in log variables to estimate the exponents (linked to Hurst index) and the other parameters.

However, since we calculate cross-covariance as well as covariances, we have an overdetermined set of equations to estimate the Hurst parameters. We experimentally show that it is preferable to eliminate this overdetermination for the estimation of the Hurst exponents. Therefore, the first conclusion of the study is that for the estimation of the Hurst exponents of the components, it is not advantageous to consider the whole multivariate process, but better to process each component separately. The second conclusion is the fact that the quality of these estimations is  almost independent of the correlation between the components. Finally we illustrate the fact that the estimation of the correlation structure is easy whereas it is very difficult to estimate the asymmetry parameters.
Our finding are based on experiments as well as theoretical proof of convergence of the estimators we exhibit. We show there almost sure convergence and provide a central limit theorem proving usual $\sqrt{n}$ convergence rate if the filters are properly chosen. 
 
 The paper is structured as follows. We present in the following section the essential facts on the mfBm needed for the paper to be self consistent. We also present the filters that we are using and the statistical properties of the filtered mfBm. Section \ref{estim:sec}
then presents the methodology we adopt. We first present basic identities highlighting the power law behavior, and then discuss the least square regression that solve our  inference problem. Section \ref{Convergence:sec} is dedicated to the theoretical study of the estimators, where we first exhibit almost sure convergence and then prove a central limit theorem. In section \ref{experiment:sec}, then, we illustrate our findings using Monte-Carlo experiments, and we present an illustration of the method on a high dimensional example. Note that the proofs of the results are given in the last section.

\section{Multivariate fBm, Filters}
\label{mfbmfact:ssec}

We recall here some basic facts about the multivariate fractional Brownian motion. For more information and proofs of the results recalled here, we refer the reader to \cite{AmblCLP10,CoeuAA10,LavaPS09,DidiP10}.

\subsection{Some facts on the mfBm}

The $p$ dimensional  multivariate fractional Brownian motion $x(t)$ is defined as a Gaussian process having stationary increments and having components jointly self-similar with parameters $(H_1,\ldots,H_p) \in (0,1)^p$. The self-similarity property can be stated as follows: for any real $\lambda > 0$, $x(\lambda t) \eqloi \lambda^H x(t)$ where $H=\mbox{diag}(H_1,\ldots,H_p) $ and $\lambda^H$ is intended in the matrix sense. The notation $\eqloi$ stands for equality of all the finite-dimensional probability distributions.

Joint self-similarity imposes many constraints on the correlation structure of the process. This has been studied in \cite{LavaPS09} where the general form of the covariance structure of a jointly self-similar process with stationary increments is obtained, without recoursing to the Gaussian assumption. This form is further studied in \cite{AmblCLP10,CoeuAA10}.  The covariance structure is shown to be characterized by $p^2$ real numbers $\rho_{ij}\in (-1,1),\eta_{ij} \in \mathbb{R},\sigma_i>0$, $i=1,\ldots,p ; j>i$. Parameter  $\sigma_i$ is  the standard deviation of the $i$th component at time 1, $\rho_{ij}$ is the correlation coefficient between the components $i$ and $j$ at time 1, and as such satisfies $\rho_{ij}=\rho_{ji}$. Parameters $\eta_{ij}$ are linked with the time-reversibility of the process. They are characterized by the antisymmetry property $\eta_{ij}=-\eta_{ji}$. 
In special cases, these parameters are known \cite{AmblCLP10,CoeuAA10}.  If the process is time-reversible, they are all equal to zero; if the process admits a causal (or an anticausal) representation, they are function of $\rho_{ij}$, $H_i$ and $H_j$. In general otherwise, they are unconstrained. 

The covariance structure of the process is as follows. The process is marginally a fractional Brownian motion. Thus 
 the covariance function of the $i$th component is the usual function \cite{MandVN68,SamoT94}
\begin{equation}
  \E [ x_i(s)x_i(t) ] \ = \frac{\sigma_i^2}{2} \left\{|s|^{2H_{i}} +
   |t|^{2H_{i}} - |t-s|^{2H_{i}}\right\}, \label{fbm:eq}
 \end{equation}
with, as mentioned, $\sigma^2_i := {\rm var}(x_i(1)) $. The cross-covariances are given by \cite{AmblCLP10} (\cite{LavaPS09} for the proof and a different parametrization)
\begin{prop} For all $(i, j)\in
\{1,\ldots,p\}^2 $, $ i\ne j$, 
\begin{eqnarray} r_{ij}(s,t) &:=& \E[  x_i(s) x_j(t) ]\\&=&
  \frac{\sigma_i\sigma_j}{2} \left\{ w_{ij}(-s) +
  w_{ij}(t)    - w_{ij}(t-s) \right\},
    \label{covariancedef:eq}
\end{eqnarray}
\end{prop}
where the function $w_{ij}(h)$ is defined  by
\begin{equation}
w_{ij}(h) = \left\{ \begin{array}{ll}
(\rho_{ij} -\eta_{ij} \mathrm{sign}(h) ) |h|^{H_i+H_j} & \mbox{ if } H_i+H_j \neq 1, \\
{\rho}_{ij} |h| + {\eta}_{ij} h \log|h| & \mbox{ if } H_i+H_j =1.
\end{array} \right.
\label{fonctionw:eq}
\end{equation}
As shown in \cite{AmblCLP10,CoeuAA10}, the form obtained for  $H_i+H_j = 1$ can be recovered by continuity from the case $H_i+H_j\not=1$. Furthermore, setting evidently $\rho_{ii}=1$ and noticing that $\eta_{ii}=0$ allows us to remark that the definition is valid if $i=j$  since it is equivalent to (\ref{fbm:eq}).

\begin{center}
\begin{figure}[p]
\includegraphics[scale=.5]{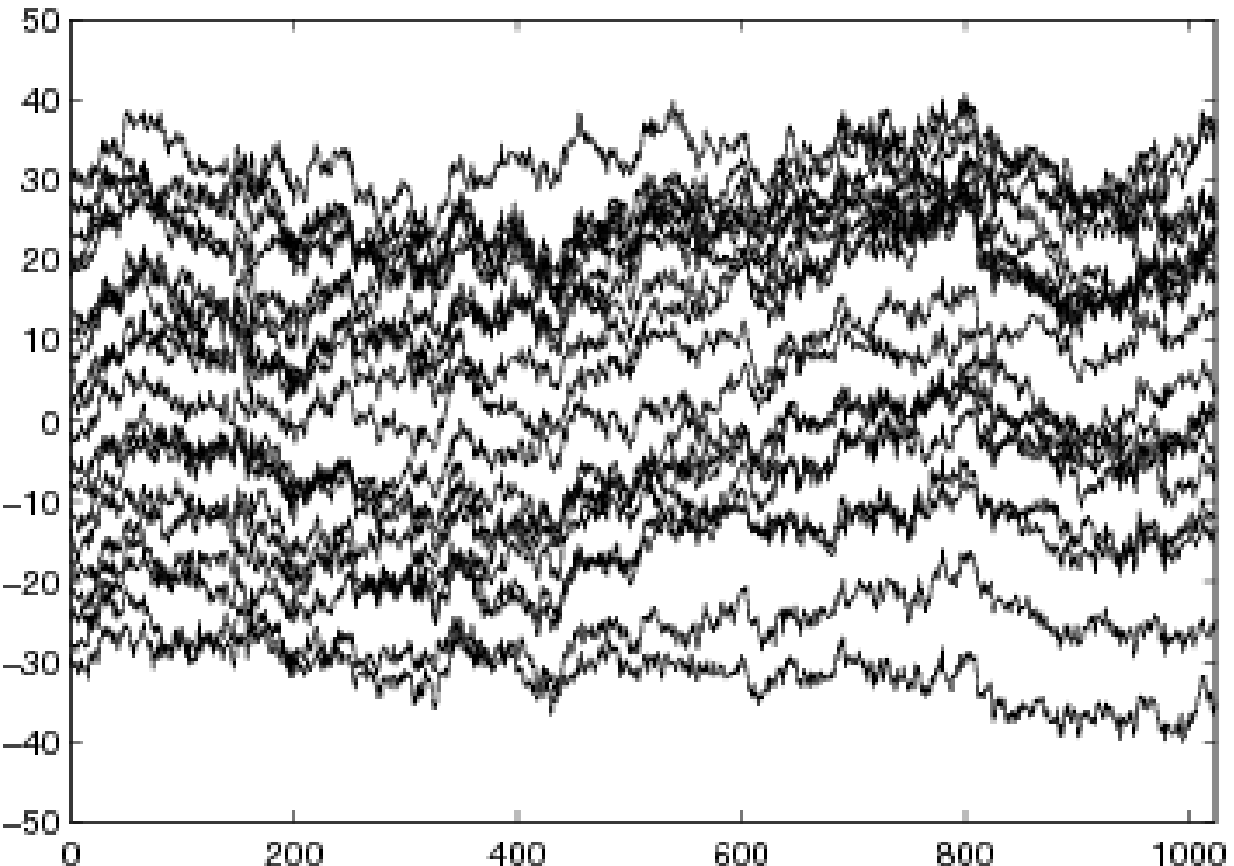} \\
\includegraphics[scale=.5]{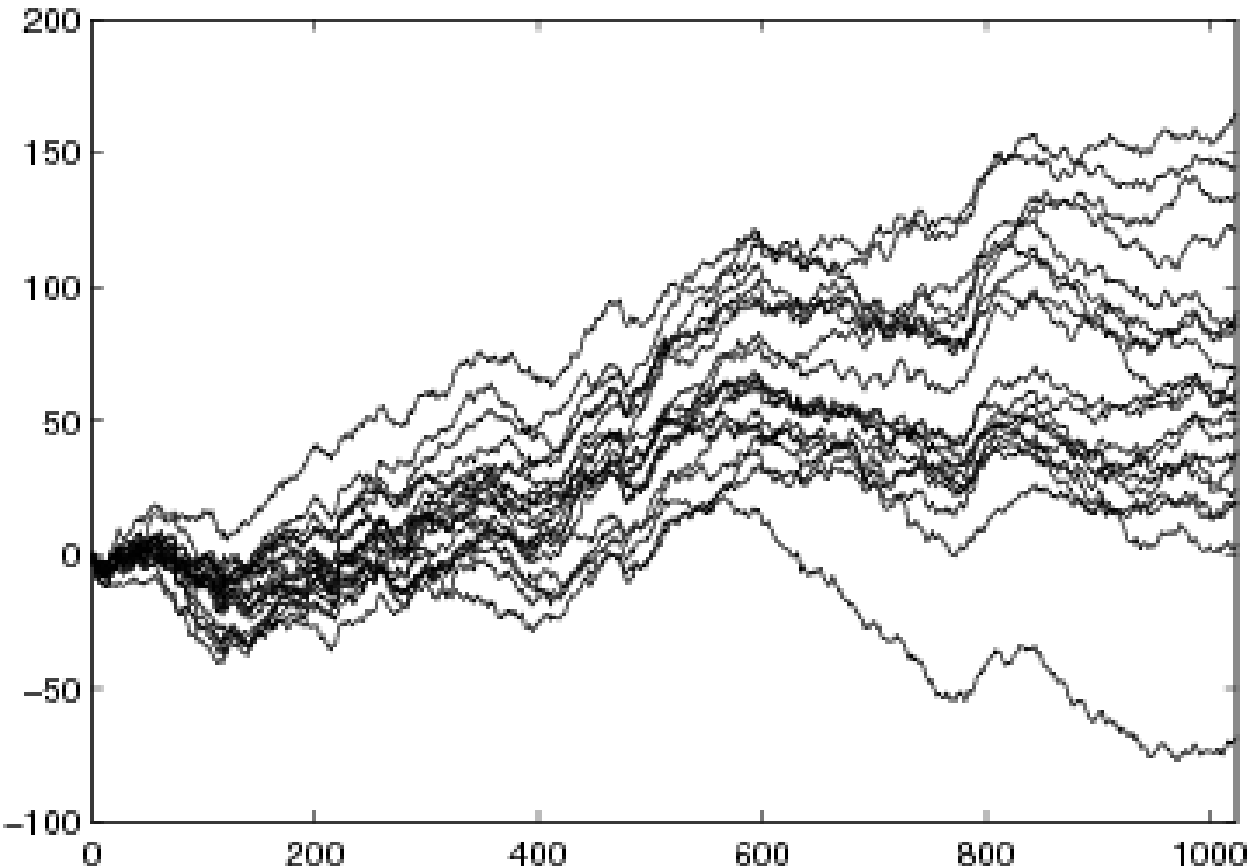}\\
\includegraphics[scale=.5]{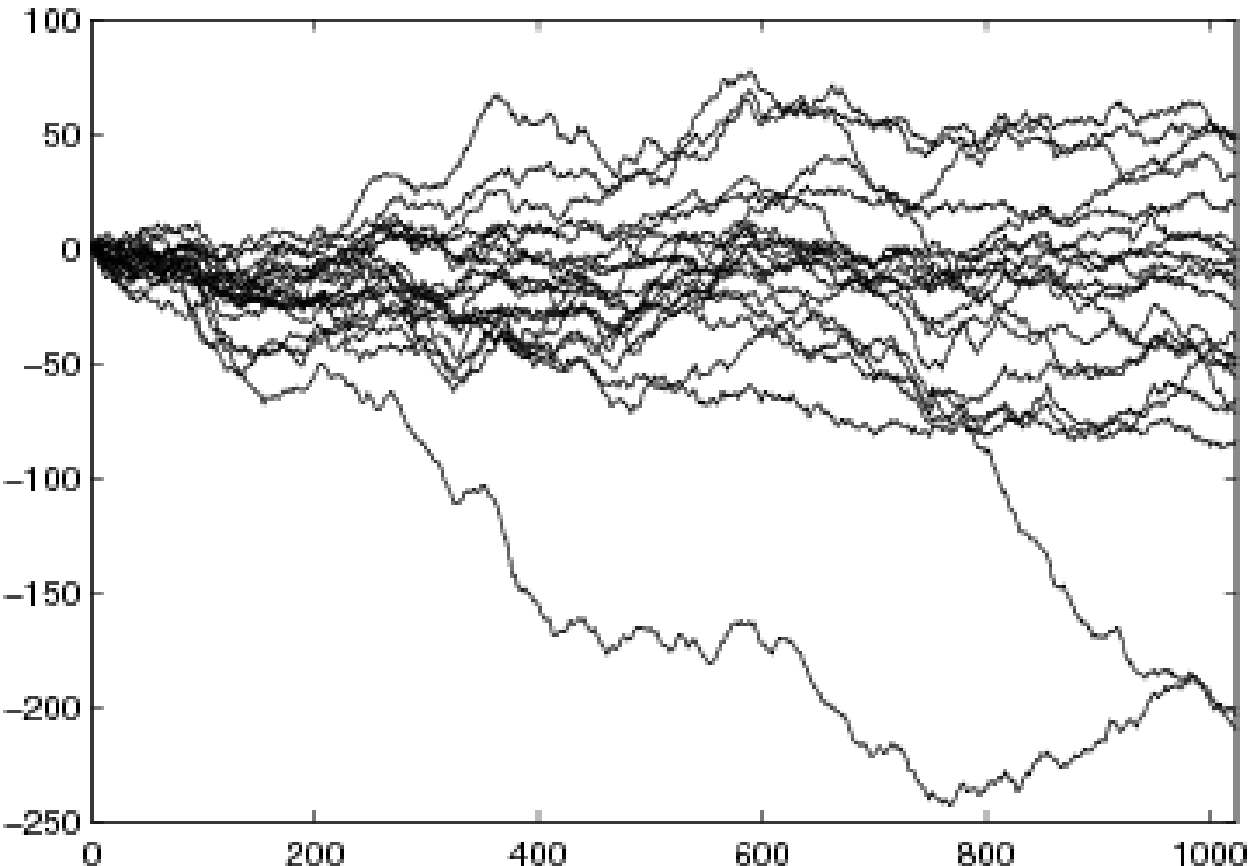}
\caption{\small Examples of discretized sample paths of a well-balanced ($\eta_{ij}=0$) mfBm of length $n=1024$, with $p=20$ components. The Hurst exponents are equally spaced in $[0.3 , 0.4]$ (upper plot), $[0.6 , 0.7]$ (middle plot) and $[0.4, 0.8]$ (bottom plot). The correlation parameters are set to 0.7 (upper and middle plot) and to 0.3 (bottom plot). The components are translated artificially in the upper plot for the sake of visibility.} 
\label{ex_mfbm:fig}
\end{figure}
\end{center}

The constraints on $\rho_{ij}$ and $\eta_{ij}$ are only necessary conditions to ensure that the matrix  given by (\ref{covariancedef:eq}) together with (\ref{fonctionw:eq}) is the cross-covariance matrix of a process. A necessary and sufficient condition has been exhibited in \cite{AmblCLP10}. This condition is the positive-definiteness of the matrix with entries
\begin{eqnarray}
\Gamma(H_i+H_j+1) \Big(\rho_{ij}\sin\big(\frac{\pi}{2}(H_i+H_j)\big) - {\bf i} \eta_{ij}\sin\big(\frac{\pi}{2}(H_i+H_j)\big), \label{condExist}
\end{eqnarray}
where ${\bf i}=\sqrt{-1}$. Interestingly, the condition  emerges when studying moving average and spectral representations of the mfBm. For example, a moving average representation can be shown to be given by (assuming $H_i\neq 1/2$ for $i=1,\ldots,p$)
\begin{equation}\label{Xma}
x_i(t)=\sum_{j=1}^p   \int_{\R} M_{i,j}^+ \left( (t-x)_+^{H_i-.5}-
  (-x)_+^{H_i-.5}\right) +
 M_{i,j}^- \left((t-x)_-^{H_i-.5} - (-x)_-^{H_i-.5}\right) W_j( \mathrm{d} x),
\end{equation}
where  $W(\mathrm{d}x) = (W_1(\mathrm{d}x), \cdots, W_p(\mathrm{d}x)) $ is a Gaussian white noise with zero mean, independent components and covariance
$\E [ W_i(\mathrm{d}x) W_j(\mathrm{d}x)] = \delta_{i,j} \mathrm{d} x $. For given parameters $\sigma, \rho, \eta$, we can find easily the terms $M_{ij}^\pm$ of the matrices $M^\pm$ (see \cite{AmblCLP10}). This representation is interesting since it shows that we have access to a whole family of different processes with different characteristics governed by the parameters. In this paper, we will for example particularly focus on the so-called causal and well-balanced cases, for which we have respectively $M^-=0$ and $M^+=M^-$. The case $M_-=0$ sets a close link between $\rho_{ij} $ and $\eta_{ij}$ whereas the case $M^+=M^-$ makes the process time-reversible leading to $\eta_{ij}=0$. The problem of simulation of such a process has been investigated in \cite{AmblCLP10} using the Chan and Wood algorithm, \cite {ChanW99}. Figure (\ref{ex_mfbm:fig}) presents some examples in order to illustrate the process.
%
%
%
%

The mfBm has by definition stationary increments. It is easy to derive the covariance structure of the increments process. Let $\Delta x(t)=x(t+1)-x(t)$ be this process (with increments of size~1) that we will refer to the multivariate fractional Gaussian noise. Then
\begin{eqnarray}
\gamma_{ij}(h) &:=& \E [ \Delta x_i(t) \Delta x_j(t+h) ] \\
&=&  \frac{\sigma_i \sigma_j}2 \bigg( w_{ij}(h-1)- 2 w_{ij}(h) + w_{ij}(h+1)\bigg).
\end{eqnarray}
The asymptotic behavior has been studied in  \cite{AmblCLP10,CoeuAA10}. We have
as $|h|\to +\infty$
\begin{equation}
\gamma_{ij}(h) \sim \sigma_i \sigma_j  |h|^{H_i+H_j-2} \kappa_{ij}(\mathrm{sign}(h))
\label{asymptotic:eq}
\end{equation}
with
\begin{equation}
\kappa_{ij}(\mathrm{sign}(h))=\left\{ 
\begin{array}{ll}
(\rho_{ij}-\eta_{ij} \mathrm{sign}(h)) (H_i+H_j)(H_i+H_j-1) & \mbox{ if } H_i+H_j\neq 1, \\ 
\tilde\eta_{ij} \mathrm{sign}(h)& \mbox{ if }H_i+H_j=1. 
\end{array}
\right.
\end{equation}
We recover here the usual behavior of the scalar fGn: each component of the mfGn can be short or long-range dependent
if its corresponding Hurst parameter is smaller or greater than 1/2, respectively. But in the multivariate case, long-range (inter)dependence can also appear in the cross-covariance. Indeed, from (\ref{asymptotic:eq}) we easily conclude that $\gamma_{ij}(h)$ is not summable as soon as $H_i+H_j\geq 1$, a case which can appear in three situations:
\begin{enumerate}
\item $H_i=1/2=H_j$
\item $H_i<1/2$ and $H_j>1-H_i$
\item $H_i>1/2$ and $H_j>1/2$
\end{enumerate}
In those cases, some troubles may appear when it comes to infer parameters of the models from data. Indeed, long-range dependence may lead to very slow convergence of estimators.

As already observed in many works \cite{Flan92,VeitA99,Coeu01}, recoursing to wavelet types of transformation is an elegant way to overcome the problem. Indeed, using wavelet types of transformation with a correctly chosen filter allows to extract the stationary part from the fBm and allows us to ``whiten" the increments. We describe such an approach in the following section.

\subsection{Discrete filtering technique and its consequence on the mfBm}

In the identification problem, we suppose to have access to a sampled version of the mfBm. We thus turn to discrete time. 
Let $\ell$ and $q$ be two positive integers.  We consider
the following set of filters ${\cal A}_{\ell,q}$:
\begin{eqnarray*}
{\cal A}_{\ell,q} = \Big\{ (a_k)_{k\in \Z}: \;   a_k=0, \; \forall k\in \Z^{-,*}\cup \{\ell+1,\ldots,+\infty\} \mbox{ and } \sum_{k\in \Z} k^l a_k=0 , \forall l =0,\ldots, q-1    \Big\}
\end{eqnarray*}
Typical examples are the difference filter $\delta_{l,0} -\delta_{l,1}$ and its compositions, Daubechies wavelet filters, and any known wavelet filter with compact support and a sufficient number of vanishing moments.

For $a\in {\cal A}_{\ell,q}$ and an integer $m\geq1$ we define the $m$th dilated version of $a$, say $a^m$ as
\begin{eqnarray*}
a^m_k = \left\{ 
\begin{array}{ll}
   a_{k/m}   &  \mbox{if } k\in m\Z  \\
    0  &    \mbox{if } k\not\in m\Z
\end{array}\right.
\end{eqnarray*}
Evidently, $a^1=a$ and $a^m \in {\cal  A}_{\ell,q}$ for any $m$. The $m$th dilated version is thus simply obtained by oversampling $a$ by a factor of $m$, {\it i.e.} by adding $m-1$ zeros between  each of the first $\ell+1$ coefficients of the impulse response $a_k$.

Let $x(t)$ be a mfBm in discrete time. We mean by this that we have at hand a collection of samples regularly taken from a continuous time mfBm. Let $x^m$ be the signal obtained by filtering $x$ with filter $a^m$. Since $x$ is multivariate, $x^m$ is also, and its components are the components of $x$ filtered by $a^m$, $x^m(t)=(x^m_1(t),\ldots, x^m_p(t))^t$ where
\begin{eqnarray*}
x^m_i(t)  = \sum_{k\in\Z} a^m_k x_i(t-k),
\end{eqnarray*}
$x$ being Gaussian with zero mean, $x^m$ is also. Now we have
\begin{eqnarray*}
\gamma^{m_1,m_2}_{ij}(h) &:=& \E [ x^{m_1}_i(t)x^{m_2}_j(t+h) ] \\
&=& \sum_{k,l\in\Z} a_k^{m_1} a_l^{m_2}   r_{ij}(t-k,t+h-l) \\
&=& -\frac{\sigma_i\sigma_j}{2} \sum_{k,l\in\Z} a_k^{m_1} a_l^{m_2} w_{ij}(h+k-l).
\end{eqnarray*}
The last equation is obtained since for any member of ${\cal A}_{\ell,q}$, $\sum_{l \in \Z} a_l=0 $. 
Using the definition of $a^m$ and of $w_{ij}$ we get
\begin{eqnarray}
\gamma^{m_1,m_2}_{ij}(h) = -\frac{\sigma_i\sigma_j}{2}  \sum_{k,l\in\Z} a_k a_l 
\big(\rho_{ij} -\eta_{ij} \mathrm{sign}(h+m_1k-m_2l) \big) \big|h+m_1k-m_2l)\big|^{H_i+H_j}
\label{covscale:eq}
\end{eqnarray}
The behavior of $\gamma^{m_1,m_2}_{ij}(h) $ has been studied in \cite{CoeuAA10} in the case of the continuous wavelet analysis of the continuous time mfBm. The result proved in Proposition 7 of the referenced paper can be developed also in the same way 
in the case of discrete wavelet transform or in the setting used here. We thus state without proof the following expansion and its consequence on the summability of $|\gamma^{m_1,m_2}_{ij}(\cdot)|^\alpha$ for $\alpha \in \N^*$.

\begin{prop}\label{prop-asympgm1m2} ${ }$\\
$(i)$ As $|h|\to +\infty$, the following equivalence holds for any $m_1,m_2\geq 1$ and any $a \in \mathcal{A}_{\ell,q}$
\begin{eqnarray*}
\gamma^{m_1,m_2}_{ij}(h) \sim  -\frac{\sigma_i \sigma_j }2 \kappa(a,q) |h|^{H_i+H_j-2q} {\tau}_{ij}(h)
\end{eqnarray*}
where $\kappa(a,q):={2q \choose q} (m_1 m_2)^{q} \left| \sum_k k^q a_k \right|^2$ and
$$
\tau_{ij}(h) = \left\{
\begin{array}{ll}
(\rho_{ij}-\eta_{ij}\sign(h)) \mbox{${ {H_i+H_j \choose 2q}}$} & \mbox{ if } i=j \mbox{ and } H_i\neq1/2 \\ 
&\mbox{ or } i\neq j \mbox{ and } H_i+H_j\neq 1 \\
\frac{\widetilde{\eta}_{ij}  \sign(h)}{2q(2q-1)} & \mbox{ if } H_i+H_j=1 \mbox{ and } H_i\neq 1/2.
\end{array} \right.
$$
$(ii)$ Let us denote by $H^\vee:=\max(H_1,\ldots,H_p)$, then for any $\alpha \in \N^*$ 
\begin{equation}\label{eq-lalpha}
q> H^\vee + \frac{1}{2\alpha} \quad \Rightarrow \quad \gamma_{ij}^{m_1,m_2}(\cdot) \in \ell^\alpha(\Z), \forall i,j=1,\ldots,p.
\end{equation}
\end{prop}

Choosing the filter $\delta_{l,0} -\delta_{l,1}$ allows us to recover~(\ref{asymptotic:eq}). The interest of filtering is revealed by taking higher order filters. Indeed, for a filter with two zero moments, the cross-covariance will be summable for all the possible values of the Hurst exponents. In some sense, the filtering aims at reducing the dependence of the cross-covariances function along time. Let us add that the key-ingredient for obtaining a central limit theorem for our proposed estimators is the square summability of all the cross-covariances functions. As stated, in~(\ref{eq-lalpha}), this will be  realized if $q=1$ and $H^\vee <3/4$ or as soon as $q\geq 2$.


We now turn to the core of the paper.

\section{Estimation method}
\label{estim:sec}

From now on, we assume having at our disposal a sample path of a mfBm (with $p> 1$ components) regularly sampled at times $t=1,\ldots,n$. For the sake of simplicity, we shall also restrict ourselves on the most interesting case $H_i+H_j \neq 1$, $\forall i,j=1,\ldots,p$.

\subsection{Basic identities}

The estimation principle relies on the covariance (\ref{covscale:eq}) for a given $m=m_1=m_2 \geq 1$. We have
\begin{eqnarray*}
\gamma^{m}_{ij}(h) := \gamma^{m,m}_{ij}(h)= -\frac{\sigma_i\sigma_j}{2}  \sum_{k,l\in\Z} a_k a_l 
\big(\rho_{ij} -\eta_{ij} \mathrm{sign}(h+m(k-l)) \big) \big|h+m(k-l))\big|^{H_i+H_j}.
\end{eqnarray*}
In particular, at lag 0 we obtain
\begin{eqnarray}
\gamma_{ii}^m (0)& =&  m^{2H_i}\sigma_i^2 \Big(-\frac{1}{2}\sum_{k,l\in\Z} a_k a_l   \big|k-l\big|^{2H_i} \Big) \label{variance:eq} \\
\gamma_{ij}^m (0)& =&m^{H_i+H_j}\rho_{ij}\sigma_i\sigma_j \Big(  -\frac{1}{2}  \sum_{k,l\in\Z} a_k a_l 
\big|k-l\big|^{H_i+H_j} \Big).   \label{correlation:eq}
\end{eqnarray}
To obtain (\ref{variance:eq}), we have made use of the fact that $\sum_{k,l=0}^\ell \sign(k-l)|k-l|^{H_i+H_j}=0$.
We note that the parameters of interest appears in the slope of the log  covariance at lag 0 when considered as a function of $\log m$. To obtain such a relation for the remaining parameters $\eta_{ij}$, we must remember that these parameters characterize the time asymmetry of the process. Since for a Gaussian process, time reversal invariance is equivalent to
$\gamma_{ij}(h)=\gamma_{ji}(h), \forall i,j$, it is tempting to extract $\eta_{ij}$ from differences like $\gamma_{ij}(h)-\gamma_{ji}(h)$.
Indeed, we have
\begin{eqnarray*}
\gamma_{ij}^m(m\ell)=m^{H_i+H_j}\sigma_i\sigma_j (\rho_{ij}-\eta_{ij} ) \Big(  -\frac{1}{2}  \sum_{k,l\in\Z} a_k a_l 
\big|\ell+k-l\big|^{H_i+H_j} \Big),
\end{eqnarray*}
where we have used the fact that the filters are zero as soon as $k>\ell$ and thus $\sign(\ell+k-l)=1$ in the double sum. 
Let us introduce the function
\begin{eqnarray*}
\pi_{ij}^a(h) :=  -\frac{1}{2}  \sum_{k,l\in\Z} a_k a_l 
\big|h+k-l\big|^{H_i+H_j},
\end{eqnarray*}
where the indices $i,j$ correspond to the fact that $\pi$ depends on the corresponding Hurst exponents. Let us underline that for all $H_i,H_j \in (0,1)$, $\pi_{ij}^a(0)>0$ for any filter $a$.
We thus have obtained the following $p^2$ equations
\begin{eqnarray}
\gamma_{ii}^m (0)& =&  m^{2H_i}\sigma_i^2 \pi_{ii}^a(0) , \hspace{.5cm}\forall i=1,\ldots,p \label{gammaii:eq} \\
\gamma_{ij}^m (0)& =&  m^{H_i+H_j} \rho_{ij}\sigma_i \sigma_j \pi_{ij}^a(0) ,  \hspace{.5cm}\forall i=1,\ldots,p, j>i \label{gammaij:eq}\\
 \gamma_{ij}^m(m\ell) -\gamma^m_{ji}(ml)&=& 2 m^{H_i+H_j}\eta_{ij}\sigma_i \sigma_j \pi_{ij}^a(\ell), \hspace{.5cm}  \forall i=1,\ldots,p, j>i \label{deltagammaij:eq} 
\end{eqnarray}

Equations (\ref{gammaii:eq}) and their wavelet counterparts in the scalar case have been used by many people to estimate the Hurst exponent ({\it e.g.} \cite{Flan92,VeitA99,Coeu01} to cite some but a few). The two others are direct extension and are going to be used in the sequel to identify the mfBm. 

At this point, the question ``which parameters do we want to estimate and how'' must be asked. If we want to only estimate Hurst parameters $H_1,\ldots,H_p$, do we have to use only the $p$ equations  (\ref{gammaii:eq}), or do we gain something by adding the $p(p-1)$ others? The parameters $H$ will be estimated by linear regression (in the log variables). Can we use these regressions to estimate the other parameters $\sigma, \rho $ and $\eta$, or is it better to consider usual empirical estimates? 

We try to adress all these questions in the following.

\subsection{Methodology}

We apply the filtering for all values of $m$ taken from a discrete set ${\cal M}$ of cardinal $|{\cal M}|$, and we thus obtain the multivariate signal  $x^m(t)$. We then evaluate the empirical estimators
\begin{eqnarray}
C^m_{ij}(h) = \frac{1}{n-m\ell -h}  \sum_{t=m\ell+ 1}^{n-h} x_i^m(t)x_j^m(t+h),
\label{def-Cijmh}
\end{eqnarray}
$C_{ii}^m(0)$ thus corresponding to the empirical moment of order~2 of the signal $x_i^m$. As we may expect that $C_{ij}^m(h)$ correctly estimates $\gamma_{ij}^m(h$), we will use this estimator to estimate the parameters of the model. Precisely, inspired by~(\ref{gammaii:eq},\ref{gammaij:eq},\ref{deltagammaij:eq}), let us introduce 
\begin{eqnarray}
\begin{array}{lcl}
    v_i^m := \log C^m_{ii}(0)   &  \mbox{ }&   \alpha_i  :=  \log \left(  \sigma_i^2 \pi_{ii}^a(0)  \right), \\
     c_{ij}^m := \log \big| C^m_{ij}(0)\big| &  & \mu_i  :=  \log \left(  \sigma_i \sigma_j |\rho_{ij}| \pi_{ij}^a(0)  \right),\\
     d_{ij}^m := \log 0.5 \big| C^m_{ij}(m\ell) -  C^m_{ji}(m\ell)\big| & &\nu_i  :=  \log \left(\sigma_i \sigma_j |\eta_{ij} \pi_{ii}^a(\ell)|  \right).
\end{array}
\label{coeffdef:eq}
\end{eqnarray}
We have to underline here that it is assumed that none of the parameters $\rho_{ij}$ and $\eta_{ij}$ is equal to zero. This could be a limitation since zero expresses the absence of correlation or the time reversibility but as we will see later the derived estimates of $\rho_{ij}$ and $\eta_{ij}$ actually do not depend on this assumption. Then, we can then write
\begin{eqnarray*}
\begin{array}{lclcl}
    v_i^m &=& 2H_i \log m  +    \alpha_i  + \varepsilon^m_{v_i} &\forall & i=1,\ldots,p  \;,\\
     c_{ij}^m &= &(H_i+H_j) \log m +  \mu_i  + \varepsilon^m_{c_{ij}} &\forall & i=1,\ldots,p ; j>i \;,\\
     d_{ij}^m &= &(H_i+H_j) \log m  + \nu_i  +\varepsilon^m_{d_{ij}}&\forall & i=1,\ldots,p ; j>i.
\end{array}
\end{eqnarray*}
The noise terms $\varepsilon$  measure the deviation of the model and can be written as
\begin{eqnarray*}
\varepsilon^m_{v_i} &=& v_i^m - \log \gamma^m_{ii}(0),\\
\varepsilon^m_{c_{ij}} &= &c_{ij}^m -\log \big| \gamma^m_{ij}(0) \big|, \\ 
\varepsilon^m_{d_{ij}}& =& d_{ij}^m  - \log 0.5 \big|  \gamma^m_{ij}(m\ell) -  \gamma^m_{ji}(m\ell)\big|.
\end{eqnarray*}

Let us now consider the vectors $H=(H_1,\ldots,H_p)^t, \alpha=(\alpha_1,\ldots,\alpha_p)^t, \mu=(\mu_{ij})^t_{i=1,\ldots,p;j>i}, \nu = (\nu_{ij})^t_{i=1,\ldots,p;j>i}$. For these two last, the ordering chosen to  create a vector is of no importance. However, to fix ideas we will use the identification $k(i,j)=(i-1)p+j-i(i+1)/2$ which corresponds to a numbering following rows. In all the following, we will often switch from matrix notation to the vector one, but the context will make it clear.

We suggest to obtain the above parameters by minimizing the following weighted mean square error objective
\begin{eqnarray*}
f(H,\alpha, \mu,\nu) = \sum_{m \in {\cal M}} \left(  w_v \sum_{i=1}^p \big(\varepsilon^m_{v_i} \big)^2
 + w_c \sum_{i=1,j>i}^p \big(\varepsilon^m_{c_{ij}} \big)^2
+w_d \sum_{i=1,j>i}^p\big(\varepsilon^m_{d_{ij}}\big)^2\right).
\end{eqnarray*}
The interest of this objective function is in the fact that it combines the three different types of ``observations", empirical variances, empirical correlation and empirical measure of asymmetry. The weights allow us to consider the advantage of including one of these types of observation in the inference problem. For example, in the usual setting, we will set $w_c=w_d=0$ and this will allow us to estimate the Hurst exponents. Considering only $w_d=0$ allows to add in the observation the empirical correlation in the hope that it will ameliorate the estimation of the Hurst exponents. 

We now solve the optimization problem 
\begin{eqnarray*}
\big(\hat{H},\hat{\alpha}, \hat{ \mu},\hat{\nu} \big) = \arg \min f(H,\alpha, \mu,\nu).
\end{eqnarray*}
The details of the calculation are provided in section (\ref{optim:ssec}).
 To write down the result, we introduce the vector of $\R^{|{\cal M}|} $,
$L:=(\log m_1,\ldots, \log m_{|{\cal M}|} )^t$. The variables $v_i,c_{ij}$ and $d_{ij}$ without exponent $m$ stand for the vectors of $\R^{|{\cal M}|} $ collecting respectively $v^m_i,c^m_{ij}$ and $d^m_{ij}$ for $ m=m_1,\ldots, m_{|{\cal M}|} )$. Furthermore, define for any vector $x \in \R^{|{\cal M}|} $ its mean $\bar{x}:=|{\cal M}|^{-1}\sum_{m \in {\cal M}}x_m$ and the centered vector $\breve{x}=x-\bar{x}$. Then, the parameters optimizing $f$ are given by 
\begin{eqnarray}
  \hat{\alpha}_k &=   &  \bar{v}_k  - 2\hat{H}_k \bar{L}, \label{alphaest:eq}\\
     \hat{\mu}_{ij} &=   &  \bar{c}_{ij}  - (\hat{H}_i+\hat{H}_j) \bar{L}, \label{muest:eq}\\
     \hat{\nu}_{ij}  &=  &  \bar{d}_{ij}  -  (\hat{H}_i+\hat{H}_j) \bar{L},   \label{nuest:eq}\\
     \hat{H}_k&=&\displaystyle \frac{\big( \breve{L}^t \breve{L}\big)^{-1}    \breve{L}^t }{\lambda }
     \Big\{   2 {v}_k  + \sum_{j\not=k} (w_c  {c}_{kj} + w_d  {d}_{kj}) \nonumber   \\ 
  & -& \displaystyle  \frac{     (w_c+w_d) 
     }{\big(\lambda+p (w_c+w_d) \big)}  \sum_{i=1}^p \Big(2 {v}_i + \sum_{j\not=i} (w_c  {c}_{ij} + w_d  {d}_{ij})\Big)\Big\}, \label{Hopt:eq}
\end{eqnarray}
where $\lambda:=4w_v+(p-2) (w_c+w_d)$. Note that setting $w_c=w_d=0$ and $w_v=1$, we find for $H_k$
\begin{eqnarray}
\hat{H}_k&=& \frac{ \breve{L}^t {v}_k}{2\breve{L}^t \breve{L}}, \label{def-Hest}
\end{eqnarray}
which is the estimator found when estimating the Hurst exponent of a scalar fBm \cite{Coeu01}. Equation~(\ref{Hopt:eq}) appears therefore as a generalization for which the estimates are still independent of the other parameters $(\sigma^2,\rho,\eta)$. When $p=2,w_v= w_c=1, w_d=0$, it takes for example the simple form:
$$
\widehat{H}_k = \frac{\breve{L}^t \big\{ 10 v_k - 2v_j + 4 c_{kj} \big\}}{ 24\; \breve{L}\breve{L}^t}, \;\; j\neq k.
$$

To conclude the identification of the mfBm, 
 parameters $\sigma_i,\rho_{ij},\eta_{ij}$, can be estimated by plugging estimators (\ref{alphaest:eq},\ref{muest:eq},\ref{nuest:eq}) into Equation~(\ref{coeffdef:eq}).
We then obtain
\begin{eqnarray}
\widehat{\sigma}_i^2 &=& \frac{e^{\hat{\alpha}_i} }{\hat{\pi}_{ii}^a(0) },   \label{def-s2iest}\\
 |\widehat{\rho}_{ij}|  &=&   \frac{e^{\hat{\mu}_{ij}} }{\widehat{\sigma}_i\widehat{\sigma}_j\hat{\pi}_{ij}^a(0)}  =  \prod_{m} \left( \frac{|C_{ij}^m(0)|}{\sqrt{C_{ii}^m(0)C_{jj}^m(0)}} \right)^{1/|\mathcal{M}|} \times \frac{\sqrt{\widehat{\pi}_{ii}^a(0) \widehat{\pi}_{jj}^a(0)} }{\widehat{\pi}_{ij}^a(0)}, \label{def-rhoijest}\\
      |\widehat{\eta}_{ij} | &=&   \frac{e^{\hat{\nu}_{ij}} }{\widehat{\sigma}_i\widehat{\sigma}_j|\hat{\pi}_{ij}^a(\ell)|}  = \frac12 \prod_m \left( \frac{|C_{ij}^m(m\ell)-C_{ji}^m(m\ell)|}{\sqrt{C_{ii}^m(0)C_{jj}^m(0)}} \right)^{1/|\mathcal{M}|} \times \frac{\sqrt{\widehat{\pi}_{ii}^a(0) \widehat{\pi}_{jj}^a(0)} }{|\widehat{\pi}_{ij}^a(\ell)|}, \label{def-etaijest}
\end{eqnarray}
where $\hat{\pi}_{ij}^a(h)$ stands for ${\pi}_{ij}^a(h)$ in which parameters $H_i$ are replaced by their estimator $\hat{H}_i$. 
or
\begin{eqnarray*}
\hat{\pi}_{ij}^a(h) =  -\frac{1}{2}  \sum_{k,l\in\Z} a_k a_l 
\big|h+k-l\big|^{\hat{H}_i+\hat{H}_j}.
\end{eqnarray*}
The r.h.s. of~(\ref{def-rhoijest}) and~(\ref{def-etaijest}) are directly obtained from~(\ref{def-s2iest}) and~(\ref{muest:eq},\ref{nuest:eq}). These last forms are very interesting, since they are not suffering from the $\log$-transformation and therefore if the model is such that $\rho_{ij}$ or $\eta_{ij}$ equals zero, such parameters can still be estimated (and thus tested). Noting from~(\ref{gammaij:eq},\ref{deltagammaij:eq}) that $\sign(\E[ C_{ij}^m(0)])=\sign(\gamma_{ij}^m(0)) = \sign(\rho_{ij})$ and $\sign(\E[ C_{ij}^m(m\ell)-C_{ji}^m(m\ell)]) =\sign( \gamma_{ij}^m(m\ell)-\gamma_{ji}^m(m\ell))=\sign(\eta_{ij})$, then estimates of $\rho_{ij}$ and $\eta_{ij}$ can be defined as follows
\begin{equation} \label{def-rhoetaijest}
\widehat{\rho}_{ij} = |\widehat{\rho}_{ij}| \times \sign(C_{ij}^m(0)) \quad \mbox{ and } 
\widehat{\eta}_{ij} = |\widehat{\eta}_{ij}| \times \sign(C_{ij}^m(m\ell)-C_{ji}^m(m\ell)),
\end{equation}
for any $m\geq 1$. Also, letting $|\mathcal{M}|=1$ allows us to recover the natural empirical estimates (obtained with one filter)
\begin{equation}\label{naturalRhoEta}
\widehat{\rho}_{ij} =\frac{C_{ij}^m(0) }{\sqrt{C_{jj}^m(0)C_{jj}^m(0)}} \times \frac{\sqrt{\hat{\pi}_{ii}^a(0)\hat{\pi}_{jj}^a(0)}}{\hat{\pi}_{ij}^a(0)} \quad \mbox{ and } \quad 
\widehat{\eta}_{ij} =\frac{C_{ij}^m(\ell) -C_{ji}^m(\ell)  }{\sqrt{C_{jj}^m(0)C_{jj}^m(0)}} \times \frac{\sqrt{\hat{\pi}_{ii}^a(0)\hat{\pi}_{jj}^a(0)}}{\hat{\pi}_{ij}^a(\ell)}.
\end{equation}
In some sense, (\ref{def-rhoijest},\ref{def-etaijest}) can be viewed, up to a factor and a sign, as the geometric mean of $\left(\frac{|C_{ij}^m(0)|}{\sqrt{C_{ii}^m(0)C_{jj}^m(0)}}\right)_{m\in \mathcal{M}}$ and $\left( \frac{|C_{ij}^m(m\ell)-C_{ji}^m(m\ell)|}{\sqrt{C_{ii}^m(0)C_{jj}^m(0)}} \right)_{m\in \mathcal{M}}$. Finally, note that we use $\hat{\pi}$ instead of $\pi$  in the renormalization of the estimators, and therefore, these estimators are coupled with the regression estimators of the Hurst exponents.

\section{Convergence analysis}
\label{Convergence:sec}
 
We concentrate in this section on the convergence of the   estimators defined by equations (\ref{Hopt:eq}), (\ref{def-s2iest}), (\ref{def-rhoijest}), (\ref{def-etaijest}) and (\ref{def-rhoetaijest}). We denote by $\theta=\left( H^t,(\sigma^2)^t,\rho^t,\eta^t\right)^t$ the whole set of parameters where the vectors $H$, $\sigma^2$,$\rho$ and $\eta$ contain all the parameters to be estimated (a vector of length $p(p+1)$), stored in an appropriate order. 
Likewise, let $\widehat{\theta}$ be the corresponding estimator with components ordered as the ones of $\theta$. Finally note our abuse of notation: the dependence on $n$, the number of observations, is not explicit in the notation of the estimators. We begin by addressing the almost sure convergence of $\widehat{\theta}$.\\

\begin{prop} \label{prop-as}
For any filter $a \in \mathcal{A}_{\ell,q}$, any set of dilations $\mathcal{M}$ and whatever the values of the weights $w_v,w_c$ and $w_d$ are, then the vector $\widehat{\theta}$  converges almost surely towards $\theta$, as $n \to  +\infty$.
\end{prop}

\bigskip
Let us underline that the almost sure convergence holds whatever the number of vanishing moments for the filter $a$ chosen, that is for all $q\geq 1$, and for all the values of the Hurst exponents. A similar result was already proved when $p=1$, {\it i.e.} for a scalar fBm (\cite{Coeu01}, Proposition 2 $(i)$). The proposition is proved in section (\ref{almostsure:ssec}). The proof relies mainly on proving
the almost sure convergence of $C_{ij}^m(h)$ to $\gamma_{ij}^m(h) $.

We now state the central limit theorem for the estimators. To state it, we need the additional definitions and notation: let
\begin{eqnarray}
C_n &:=& \left(
C_{ii}^m(0)  , C_{ij}^m(0), C_{ij}^m(m\ell), C_{ij}^m(-m\ell), \; (m\in \mathcal M,i,j=1,\ldots,p, j>i)
 \right)^t \label{def-Cn}\\
\gamma=\E[C_n] &=& \left(\gamma_{ii}^m(0)  , \gamma_{ij}^m(0), \gamma_{ij}^m(m\ell), \gamma_{ij}^m(-m\ell), \; (m\in \mathcal M,i,j=1,\ldots,p, j>i) \right)^t. \nonumber
\end{eqnarray}
In these notation, the vectors of length $D_{\mathcal{M},p}:=|\mathcal{M}|p+3|\mathcal{M}|p(p-1)/2=|\mathcal{M}|p(3p-1)/2$ are ordered as follows: first we put the empirical variances starting with the component $x_1$ for all the values of $m$ and then the component $x_2$,\ldots. Then, we place the empirical covariances at lag zero (with the same convention: $i,j$ fixed, and $m$ is varying), then the empirical covariances at lag $m\ell$ and finally the ones at lag $-m\ell$.

Now, we note that the definition of the estimators by equations (\ref{Hopt:eq},  \ref{def-s2iest}, \ref{def-rhoijest}, \ref{def-etaijest}, \ref{def-rhoetaijest}) define unambiguously a function  $g: \R^{D_{\mathcal{M},p}} \rightarrow \R^{ p(p+1) }$ that maps the vector $C_n$ to the vector $\widehat{\theta}$. In the following result, the notation $\stackrel{d}{\to}$ stands for the convergence in distribution as $n\to +\infty$. Recall that the quantity $H^\vee$ denotes the largest Hurst exponent, that is $H^\vee:=\max_{i=1,\ldots,p} H_i$.\\

\begin{prop} Under the notation and assumptions of Proposition~\ref{prop-as} with the order of the filter, $q$, satisfying $q>H^\vee+1/4$, then\\
$(i)$ \begin{equation}\label{normCn}
\sqrt{n}(C_n-\gamma)\stackrel{d}{\to} \mathcal{N}(0,\Sigma),
\end{equation}
where $\Sigma$ is the $D_{\mathcal{M},p} \times D_{\mathcal{M},p}$ matrix explicitly given by~(\ref{eq-Sigma}) and (\ref{matrice-cov-def:eq}) p.\pageref{matrice-cov-def:eq}.\\
$(ii)$ The vector $\theta$ satisfies $\theta=g(\gamma)$, $g$ is differentiable in $\gamma$ and 
\begin{equation} \label{normTheta}
\sqrt{n} \left( \widehat{\theta} - \theta \right)
\stackrel{d}{\longrightarrow} {\cal N} \Big(0, \nabla g(\gamma)\Sigma  \nabla g(\gamma)^t\Big),
\end{equation}
where $\nabla g(\gamma)$ denotes the gradient of $g$ (a $p(p+1) \times D_{\mathcal{M},p}$ matrix) evaluated at $\gamma$.  
\label{prop-norm}
\end{prop}

\bigskip

The matrix $\Sigma$ is quite complex but it may be evaluated (or at least approximated because it containes infinite series) in order to build asymptotic confidence intervals. The most interesting point of this result is that the rate of convergence of $\widehat{\theta}$ is the optimal one, $\sqrt{n}$ (for the whole set of parameters and whatever the values of the weights $w_v,w_c,w_d$) as soon as $q\geq 2$ (see also remark after Proposition~\ref{prop-asympgm1m2}). The proof of this proposition is given in section  \ref{clt-proof:ssec}. It mainly consists in establishing $(i)$. The second point $(ii)$ will be derived using the classical delta method. We underline that no assumption is made on $\gamma$ in $(ii)$. This means that the differentiability is true for all the values of the parameter vector $\theta$ (such that the models exists). In particular, there is no differentiability problem when $\rho_{ij}$ and/or $\eta_{ij}$ equals zero. This may allow us to use~(\ref{normTheta}) to test the absence of correlation at lag zero or to test the time reversibility of the process.

In the next section, we show via experiments that when interested in estimating solely the Hurst exponents, the best strategy is to set $w_c$ and $w_d$ to 0. In this case, the gradient vector of $g$ evaluated at $\gamma$ reduces to terms of the form $\frac{\breve{L}^t}{2\breve{L}^t\breve{L}} \frac{1}{\gamma_{ii}^m(0)}$ and, with little algebra, we may derive the more simple and nice central limit theorem
\begin{eqnarray}
\sqrt{n}(\widehat{H}-H) \stackrel{d}{\longrightarrow}  {\cal N} \left(0, \frac{1}{4 (\breve{L}^t\breve{L})^2 } (I_p\otimes \breve{L})^t \widetilde\Sigma (I_p\otimes \breve{L}) \right)  \label{TCLHest}
\end{eqnarray}
where $\widetilde\Sigma$ is the  $Mp\times Mp$ matrix with elements
\begin{eqnarray*}
 2 \sum_{j\in\Z} \frac{ {\gamma}_{i_1i_2}^{m_{1},m_{2}}(j)^2 }{\gamma_{i_1i_2}^{m_{1}}(0)\gamma_{i_1i_2}^{m_{2}}(0)}.
\end{eqnarray*}
And we end this section by remarking that when $p=1$, (\ref{TCLHest}) is in agreement with Proposition~4 (ii) obtained in~\cite{Coeu01}.


\section{Experiments} \label{experiment:sec}

The previous theoretical results are important since they insure convergence at a good rate. However, the complexity of the variance terms makes the results difficult to exploit, especially when we want to compare different estimators corresponding to different choices of the weights $w_{v,c,d}$. Hence, we now turn to some Monte-Carlo experiments to study the estimators, and we present an illustration of the method on a high dimensional example.

\subsection{Experimental study of convergence}
Depending on the weights used in the objective function, we may study 8 different situations. However, all of them are not  {\it a priori} useful and do not perform well ({\it a posteriori}). For example, setting $w_v=0$ leads to very poor performance and the corresponding cases are not studied. 
We limit ourselves to the  three situations for the estimation of the Hurst exponents:
\begin{enumerate}
\item $w_v=1$, $w_c=w_d=0$: this case corresponds to applying the univariate estimator to each of the components of the mfBm.  To indicate this situation, we add a $v$ in exponent to the estimators. For example we will study $\hat{H}^v$.
\item  $w_v=1=w_c$, $w_d=0$: in this case, we study the advantage of including the empirical correlation in the observation set of the linear regression. To indicate this situation, we add a $c$ in exponent to the estimators. 
\item $w_v=1=w_c=w_d$: in this setting, the empirical measure of asymmetry is also taken into account in multiple regression. To indicate this situation, we add a $d$ in exponent to the estimators. 
 \end{enumerate}

We have generated 100 snapshots of length $n=1000$ for different cases of mfBm: different Hurst exponents, different dimensions, different correlation coefficients, for the causal, well-balanced and general mfBm. For the causal case, parameter $\eta_{ij}$ is dependent on $\rho_{ij}$, in the well-balanced case it is zero, and we have set it to $0.2\times (1-H_i-H_j)$ for all $i,j$ in the general case. The range of the parameters is constrained by the existence conditions recalled by equation (\ref{condExist}).
For each case described, we have evaluated the Mean Square Error (MSE) of the estimators in the $v,d,c$ cases reported above. For the correlation and the asymmetry coefficient, we have studied the $v$ and $d$ cases only, as well as the empirical estimators, in which parameters $\hat{H}^v$ are used in the renormalization. The hyperparameters used in this simulation study are $\mathcal{M}=\{1,\ldots,5\}$ and the generic filter $a=db4$ corresponding to a wavelet Daubechies filter with two zero moments. These choices are guided by the fact that they provided good results when dealing with a monovariate fBm, \cite{Coeu01}. Other parameters have been tried leading to the same general conclusions.

\begin{center}
\begin{figure}[p]
\includegraphics[scale=.55]{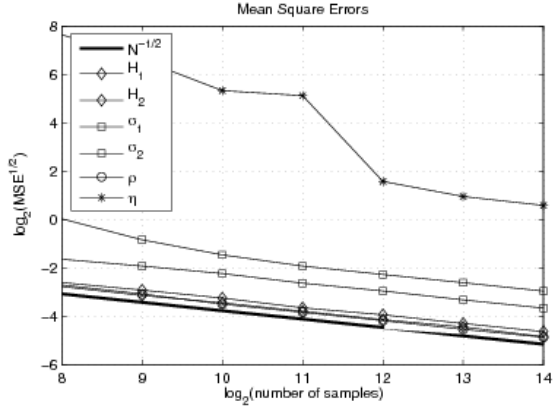}
\caption{\small Estimated standard deviation of the estimators obtained from the regression as a function of the size of the sample. The plot is a log-log plot. The estimators of $H$ is $\hat{H}^v$, and this is used in the estimation of the  others. We observe a clear $1/\sqrt{n}$ behavior. This rate is confirmed by the theoretical analysis. Note that this result is not clear for $\hat{\eta}$. The sample size should be much gretear to validate or invalidate this rate. Parameters chosen here : $H_1=0.3,H_2=0.8,\sigma_1=2,\sigma_2=1,\rho_{12}=0.4,\eta=0$.} 
\label{converge:fig}
\end{figure}
\end{center}

The results are reported in tables \ref{hest:tb},\ref{rhoest:tb} and \ref{etaest:tb}. The main conclusions from these experiments are the following:
\begin{itemize}
  \item Regarding the estimation of the Hurst exponents, adding the empirical correlation as an observation over which regression is performed does not improve the performance of the scalar estimator applied to each of the components. If $\hat{H}^v$ and $\hat{H}^c$ perform equally well   for high correlation coefficient, the performance of the latter is at least one order of magnitude less than the performance of the former when the correlation coefficient goes to zero. The performance of $\hat{H}^v$ appears almost independent of $\rho$. The same conclusion holds for $\hat{H}^d$. However, the estimators including the empirical asymmetry is considerably degraded, at least two orders of magnitude worse. 
  \item The estimation of $\rho_{ij}$ using the renormalized empirical estimator or the regression estimator based on the variance data only (when plugging in the estimate of $H$) leads to the same level of performance. 
  \item Parameters $\eta_{ij}$ are very difficult to estimate, at least with the method adopted here. The difficulty is conformed in figure \ref{converge:fig} where the MSE is plotted in a log-log plot of the estimated standard deviation versus the sample size. If 
  $\sqrt{n}$ is clearly observed for the other estimators, it is not (almost) observed by $\hat{\eta}_{ij}$, at least for the sample size up to $2^{14}$ points. 
\end{itemize}

Thus, in the following, we focus on the convergence of the estimators of the Hurst exponents, the correlation and the variance
$(\hat{H} , \sigma_i,\rho_{ij},\eta_{ij})^v$.

\begin{center}
\begin{figure}[p]
\includegraphics[scale=.4]{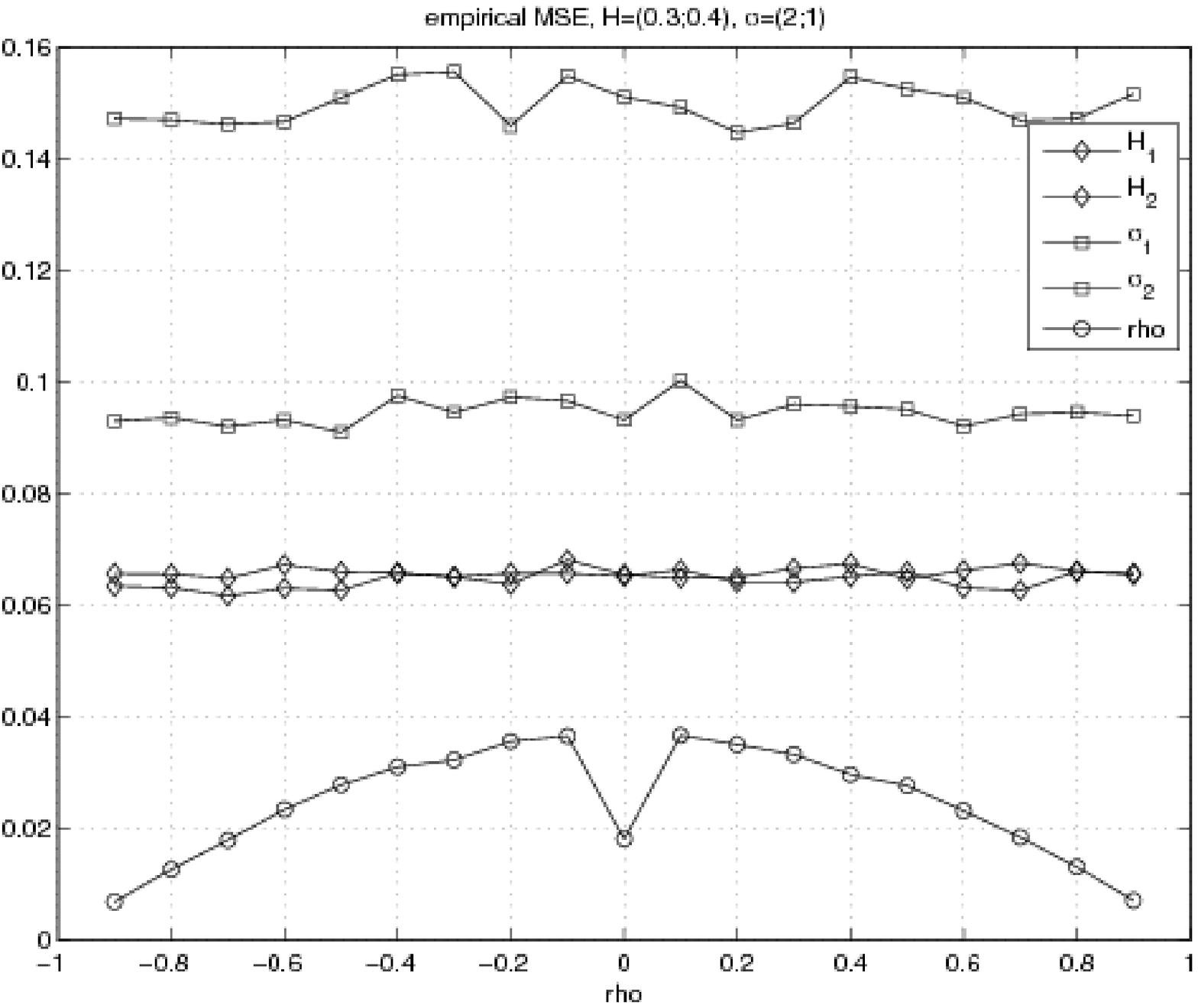} \includegraphics[scale=.4]{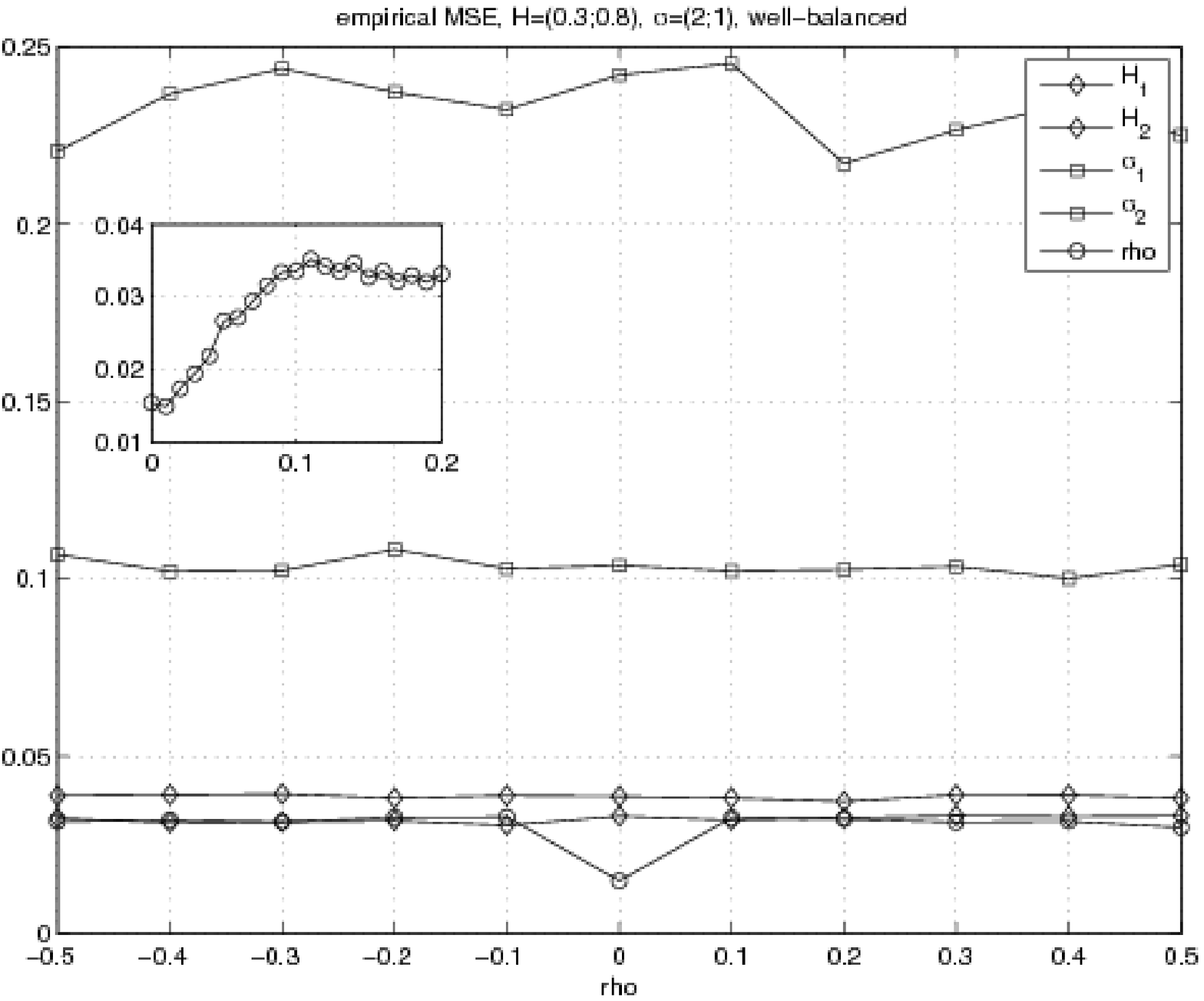}
\caption{\small Estimated mean square error of the estimators obtained from the regression as a function the correlation coefficient. Left plot: $H_1=0.3,H_2=0.4,\sigma_1=2,\sigma_2=1$, causal case. Right plot: $H_1=0.3,H_2=0.8,\sigma_1=2,\sigma_2=1$ well-balanced case. The inset depicts a zoom for the MSE of $\hat{\rho}$ {\it vs} $\rho$. Note that in the right plot, $\rho$ can not vary in the whole interval $(-1,1)$ because otherwise the model is undefined. } 
\label{msevsrho:fig}
\end{figure}
\end{center}

We have already remarked that the estimators of the Hurst exponents seem almost independent of the correlation. We thus study the behavior of the estimators with respect to the correlation. For $p=2$, we use Monte-Carlo simulation (1000 snapshot of 1024 samples each, $m=5$ dilations used) to plot the MSE of each estimator as a function of $\rho_{12}$. Results are displayed in figures 
(\ref{msevsrho:fig}). In the left plot, we study the causal case for  $H_1=0.3$ and $H_2=0.4$, whereas the right plot is concerned with the well-balanced case with  $H_1=0.3$ and $H_2=0.8$.  In the left plot, the admissible range of $\rho$ is almost all the interval $(-1,1)$, whereas it is restricted to approximately $(-0.5,0.5)$ in the case of the right plot.

The main conclusion of these plots is the fact that the estimation of the Hurst exponents and of the variances are almost not dependent on the correlation coefficient between the components of the mfBm. However, the quality of the estimation of the correlation coefficient depends on the actual value of the coefficient, and depends on it in a rather strange non monotone way. Indeed, the MSE increases with $\rho$ for moderate values and then decreases.

\subsection{A high dimensional example} \label{high}

As a conclusion we present an illustration of the model and its identification in the context of complex networks. Suppose that we observe a $p=100$ dimensional fBm obtained from a graph as follows. Let $A$ be the lower triangular part of the adjacency matrix of the graph ({\it i.e.} $A_{ij}=1 \Leftrightarrow$ nodes $i$ and $j<i$ are connected, and $A_{ij}=0$ otherwise). Then, each non zero elements is given a random value, and we use $\rho= (I-A)^{-1}(I-A)^{-t}$ as
the correlation matrix of the mfBM. The rationale hidden there is the following. Let $X$ be a 100 dimensional vector such that $X=AX+B$ where $B$ is 100 dimensional Gaussian random vector with zero mean and identity correlation matrix. Then $X$ has obviously $\rho$ as correlation matrix. 

\begin{center}
\begin{figure}[p]
\includegraphics[scale=.33]{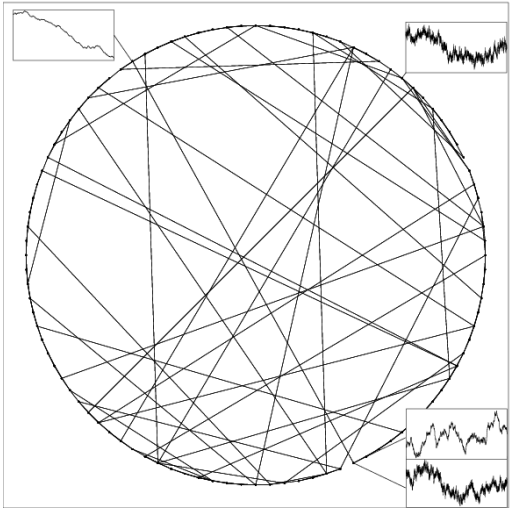}
\caption{\small Watts-Strogatz   graph used to model correlation between the components of a 100 dimensional mfBm. The two south-east inset depict components 18 and 19 that are in the example negatively correlated $\rho_{18,19}=-0.09$, with $H_{18}=0.64$ and $H_{19}=0.30$. The north-west inset depict component 76, positively correlated with 19, $\rho_{76,19}=0.13$, with Hurst exponent $H_{76}=0.7$.} 
\label{strogfBm:fig}
\end{figure}
\end{center}

As underlying graph, we choose a Watts-Strogatz model.
 A Watts-Strogatz network is a model of complex network that jointly presents the property of small-world effect (small mean geodesic distance) and the property of high clustering (neighbours of a node are strongly connected) \cite{WattS98}. This model was one of the first that adequately described graphs with these two properties. It is in a sense in between Erd\"os-R\'enyi random graph (low clustering and small mean geodesic distance) and regular grids (high clustering but low small mean geodesic distance). It is obtained by randomly rewiring edges in a regular grid. In the example depicted here (see figure \ref{strogfBm:fig}), we use a ring of nodes where each nodes is connected backward and forward with two neighbors, and each edges is rewired to a randomly chosen node (possibly the same)  independently of the others with probability 0.2 (self-connections are prohibited).
This gives the adjacency matrix we use to create a $100$ dimensional fBm as described above. 

The resulting sample paths are illustrated in figure 
(\ref{strogfBm:fig}) where we have plotted in some insets some components. For example, component 19 with Hurst exponent 0.3 is positively correlated with component 76, which Hurst exponent is slightly greater than  0.7, but negatively correlated with 18 which Hurst exponent is 0.64. 
Furthemore, since $H_{19}+H_{76}>1$, the two components are long-range cross-correlated, whereas 18 and 19 are short-range correlated.

\begin{center}
\begin{figure}[p]
\includegraphics[scale=.55]{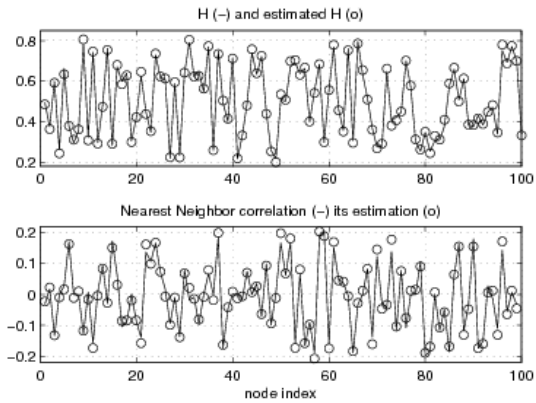}
\caption{\small Estimation of the Hurst exponents and of $\rho_{i,i+1}$ for the high dimensional example described in Section~\ref{high}.} 
\label{EstHCor:fig}
\end{figure}
\end{center}

\begin{center}
\begin{figure}[p]
\includegraphics[scale=.55]{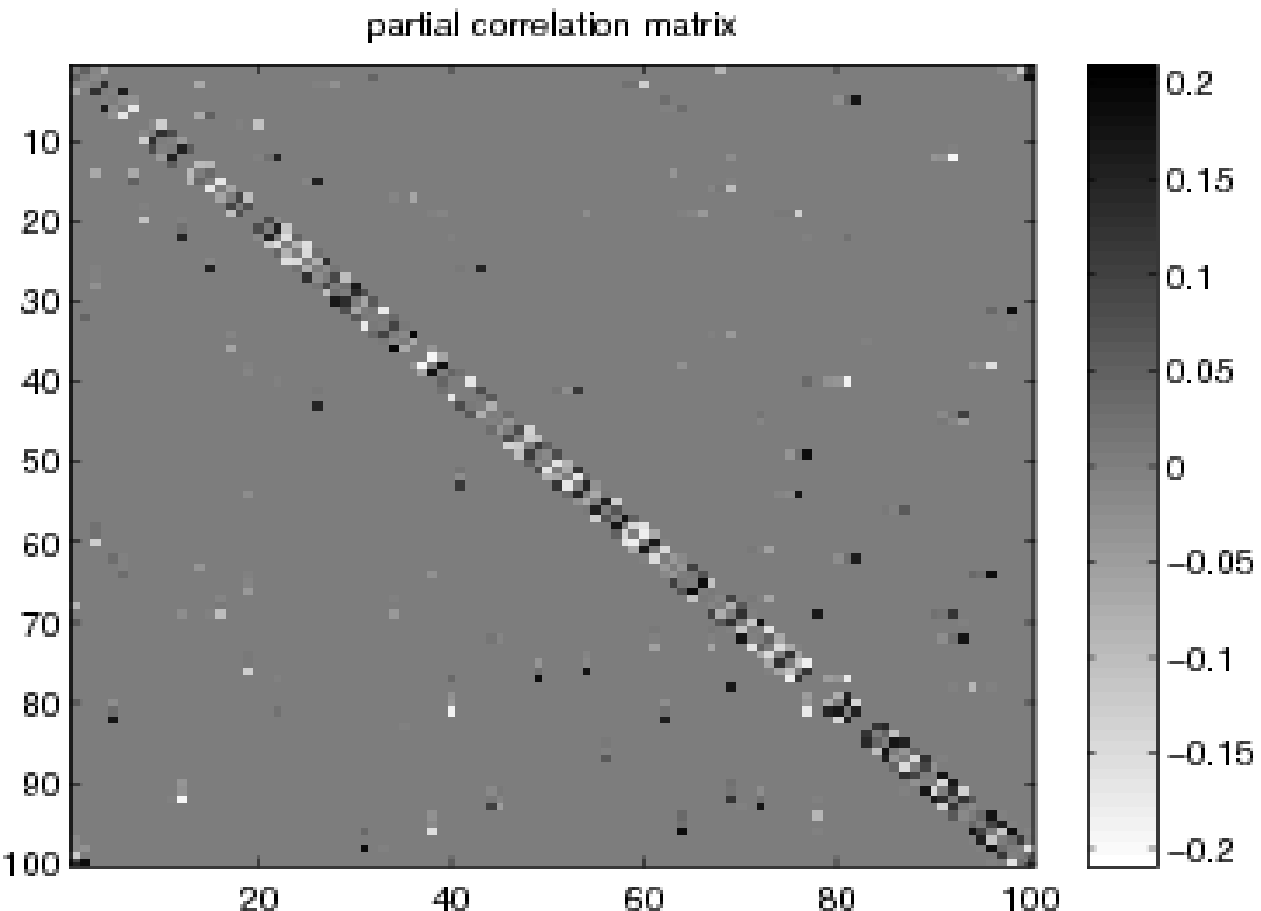} \\
\includegraphics[scale=.55]{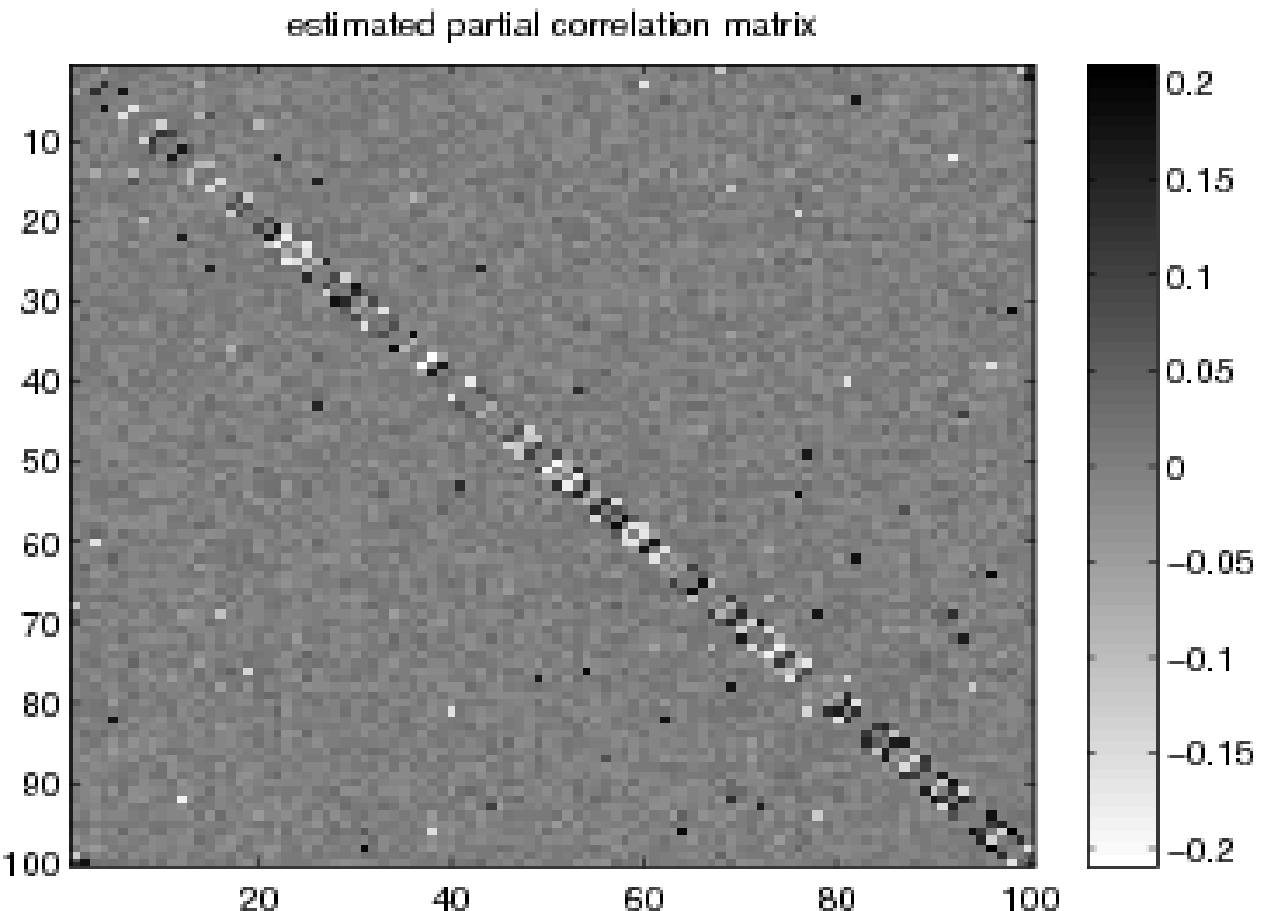}
\caption{\small True and estimated partial correlation matrix as the inverse of the correlation matrix for the high dimensional example described in Section~\ref{high}.} 
\label{ParCor:fig}
\end{figure}
\end{center}

We have generated a sample path of length 8192 samples, on which we apply our estimation procedure. Since the procedure does not depend on the dimension, we of course obtain good results for the estimation of the parameters. In figure (\ref{EstHCor:fig}), we plot the true vector $H$ and its estimation, as well as the correlation $\rho_{i,i+1}$ and its estimation. We also show in figure (\ref{ParCor:fig}) the true partial correlation matrix and its estimation {\it via }  the inverse of $\widehat{\rho}$. Recall that in the Gaussian case,  a zero partial correlation between two components is equivalent to the independence between the two components conditionnally to the remaining components.  This could be used to infer dependence link between the components of the process, as is done for example in \cite{AchaSWSB06} for connectivity studies in the brain \cite{Spor10}.

\section{Proofs} \label{sec-proofs}

\subsection{Optimization of $f$}
\label{optim:ssec}

We differentiate $f$ with respect to all the parameters and set the derivatives to zero to obtain necessary conditions for optimality. For parameters $\alpha$ we solve and obtain immediately
\begin{eqnarray}
\hat{\alpha}_k &=& \frac{1}{|{\cal M}|} \sum_{m \in {\cal M}}  [v_k^m -2\hat{H}_k \log m] \\
&=&  \bar{v}_k  - 2 \hat{H}_k \bar{L},  \label{alphaestim:eq}
\end{eqnarray}
where $L=(\log m_1,\ldots, \log m_{|{\cal M}|} )^t$, $v_k = (v_k^{m_1},\ldots,v_k^{m_{|{\cal M}|}}  )^t$,  and $\bar{x}_j=|{\cal M}|^{-1}\sum_{m \in {\cal M}}x_j^m$. Likewise, we easily get
\begin{eqnarray}
    \hat{\mu}_{ij} &=   &  \bar{c}_{ij}  - (\hat{H}_i+\hat{H}_j) \bar{L}\label{muestim:eq} \\
     \hat{\nu}_{ij}  &=  &  \bar{d}_{ij}  -  (\hat{H}_i+\hat{H}_j) \bar{L},  \label{nuestim:eq}
\end{eqnarray}
where $\mu_{ij} = (\mu_{ij} ^{m_1},\ldots,\mu_{ij} ^{m_{|{\cal M}|}}  )^t$.
Obtaining parameters $\hat{H}_k$ requires a little bit more work. Differentiating $f$ with respect to $H_k$ and setting the result to zero leads
to
\begin{eqnarray*}
\frac{1}{|{\cal M}|} \sum_{m \in {\cal M}}  w_v (v_k^m -2 H_k \log m -\alpha_k)(2\log m)  &+& \\
\frac{1}{|{\cal M}|} \sum_{m \in {\cal M}} w_c\sum_{j\not=k} \big( c_{kj}^m -(H_j+H_k) \log m -\mu_{kj})\log m &+&\\
\frac{1}{|{\cal M}|} \sum_{m \in {\cal M}} w_d\sum_{j\not=k} \big( d_{kj}^m -(H_j+H_k) \log m -\nu_{kj})\log m &=& 0.
 \end{eqnarray*}
Introducing the centered vector $\breve{x}=x-\bar{x}$ and replacing in the previous equation parameters $\alpha_k$, $\mu_{kj}$ and $\nu_{kj}$ by their estimate (\ref{alphaestim:eq},\ref{muestim:eq},\ref{nuestim:eq}), we obtain 
$$2w_v \big(  \breve{L}^t \breve{v}_k -2H_k \breve{L}^t\breve{L} \big) +
w_c \big(  \breve{L}^t \breve{\mu}_{ij} -(H_j+H_k) \breve{L}^t\breve{L} \big) +
w_d \big(  \breve{L}^t \breve{\nu}_{ij} -(H_j+H_k) \breve{L}^t\breve{L} \big) =0.
$$
Isolating $H$ terms leads to
\begin{eqnarray*}
H_k ( 4w _v+(p-2) (w_c+w_d) ) +(w_c+w_d)  \sum_{j=1}^d H_j &=&\big( \breve{L}^t \breve{L}\big)^{-1}    \breve{L}^t   \Big\{   2 \breve{v}_k  + \sum_{j\not=k} (w_c  \breve{c}_{kj} + w_d  \breve{d}_{kj}) \Big\}
\end{eqnarray*}
Collecting the $p$ equations into a vector, we have to solve
\begin{eqnarray}
\big( ( 4w _v+(p-2) (w_c+w_d) ) I_p +(w_c+w_d)J_p \big) H = X   \label{tosolve:eq}
 \end{eqnarray}
where the $k$th component of $X$ is $ \big(\breve{L}^t \breve{L}\big)^{-1}  \breve{L}^t   \Big\{   2 \breve{v}_k  + \sum_{j\not=k} (w_c  \breve{c}_{kj} + w_d  \breve{d}_{kj}) \Big\}$, where $I_p$ is the $p$ dimensional identity matrix and $J_p$ is the $p$ dimensional matrix which entries are all equal to 1. Now consider the auxiliary result
 \begin{lem} \label{lem-aux}
Let $\lambda,\lambda^\prime>0$. The $(p,p)$ matrix $B=\lambda I_p + \frac1p (\lambda'-\lambda)J_p$ has eigenvalues $\lambda$ and $\lambda'$  of respective multiplicity $p-1$ et $1$. The inverse of $B$ is thus
$$
B^{-1} = \frac1\lambda I_p + \frac1p \left( \frac1{\lambda'} -\frac1\lambda \right) J_p.
$$
\end{lem}
Applying Lemma~\ref{lem-aux} (for which the proof is omitted) to (\ref{tosolve:eq}) and noticing that for any $z\in \R^{|\mathcal{M}|}$, $\breve{L}^t \breve{z}=\breve{L}^t z$, we obtain~(\ref{Hopt:eq}).

\subsection{Proof of Proposition~\ref{prop-as}}
\label{almostsure:ssec}

\begin{proof}
The only thing to prove is that for fixed $h$, for all $i,j=1,\ldots,p$ and $m \in \mathcal{M}$,
\begin{equation}\label{convCijmh}
C_{ij}^m(h) \stackrel{a.s.}{\to} \gamma_{ij}^m(h),
\end{equation}
as $n\to +\infty$ (the notation $\stackrel{a.s.}{\to}$ stands here and in the following for the almost sure convergence). Indeed, if (\ref{convCijmh}) is true, the following convergences hold
\begin{eqnarray*}
v_i^m &\stackrel{a.s.}{\to} & \log \gamma_{ij}^m(0) = 2H_i \log m + \alpha_i, \\
c_{ij}^m &\stackrel{a.s.}{\to} & \log |\gamma_{ij}^m(0)| = (H_i+H_j)\log m + \mu_{ij}, \\
d_{ij}^m &\stackrel{a.s.}{\to} & \log 0.5|\gamma_{ij}^m(m\ell)-\gamma_{ji}^m(m\ell)| = (H_i+H_j)\log m + \nu_{ij}. 
\end{eqnarray*}
Plugging these results in~(\ref{Hopt:eq}) and noting that $\breve{L}^t 1=0$ leads with a little computation to the convergence of $\widehat{H}_k$ to $H_k$ for all $k$ and then to the convergence of $\widehat{\alpha}_k$ to $\alpha_k$. The convergences of $\sigma^2,\rho$ and $\eta$ follow from their respective definition and~(\ref{convCijmh}) applied to $h=0,\pm m\ell$. 

Let us now focus on the proof of~(\ref{convCijmh}). Define $y(k)=x_i^m(k)x_j^m(k+h)$ and assume that $y(\cdot)$ is observed at times $1,\ldots,n$. This is not a loss of generality since for fixed $m,\ell$ and $h$,  $n-m\ell-h \sim n$ as $n\to +\infty$. Let $y:=E[Y(k)]=\gamma_{ij}^m(h)$ and $\overline{y}_n:= n^{-1} \sum_{k=1}^n y(k)-y$. From Theorem 6.2 of \cite{Doob53}, p. 492, establishing a condition under which almost sure convergence is implied by mean-squared convergence for the convergence of empirical means of discrete stationary processes, the proof will be ended if we manage to prove that $E[\overline{y}_n^2]=o(1)$.  
Since $y(\cdot)$ is a stationary sequence
\begin{equation}\label{varybn}
E[\overline{y}_n^2] = \frac1{n^2}\sum_{k=1}^n (n-1-|\tau|) r_y(\tau) \leq \frac1n \sum_{\tau=1}^n  |r_y(\tau)|,
\end{equation}
where $r_y$ is the covariance function of $y(\cdot)$ given by $r_y(\tau):=E[y(k)y(k+\tau)]$. Using for example Isserlis formula, \cite{Iss18}, we can derive
\begin{equation}\label{rytau}
r_y(\tau) = \gamma_{ii}^m(\tau) \gamma_{jj}^m(\tau) + \gamma_{ij}^m(\tau+h) \gamma_{ji}^m(\tau-h).
\end{equation}
Proposition~\ref{prop-asympgm1m2} $(i)$ states in particular that $\gamma_{ij}(\tau) = \mathcal{O}(|\tau|^{H_i+H_j-2q})$ as $|\tau|\to +\infty$. Moreover, let us recall that for $H\in (0,1)$
\begin{equation} \label{suite}
\frac1n \sum_{|\tau|\leq n} \frac1{(1+|\tau|)^{2(2H-2q)}} = \left\{ \begin{array}{ll}
\mathcal{O}(1/n) & \mbox{ if } q\geq2 \mbox{ or } q=1 \mbox{ and } H<3/4. \\
\mathcal{O}(\log n/n) & \mbox{ if } q=1 \mbox{ and } H=3/4. \\
\mathcal{O}(1/n^{2-2H}) & \mbox{ if } q=1 \mbox{ and } H>3/4.
\end{array} \right.
\end{equation}
From (\ref{varybn}) and (\ref{rytau}), then using Cauchy-Schwartz inequality and (\ref{suite}) (with $H=H_i,H_j,(H_i+H_j)/2$) allows us to conclude that $E[\overline{y}_n^2]=o(1)$.
\end{proof}

\subsection{Proof of Proposition~\ref{prop-norm}}
 \label {clt-proof:ssec}
 
 The proof of the central limit theorem is done in three steps. First we prove in Lemma \ref{lem-normCijm} a central limit theorem for $C_{ij}^m(h)$, then a central limit theorem for the vectors containing all the data used in the regression, {\it i.e.}  $C_{ij}^m(0)$, $C_{ij}^m(ml)$ and $C_{ji}^m(ml)$, i.e. the proof of Proposition~\ref{prop-norm} $(i)$. Finally we apply the delta method to prove the central limit theorem for the estimators. 
 We thus begin with the first limit theorem, stated as a lemma:
\begin{lem} \label{lem-normCijm}
Let $i,j=1\ldots,p$, $m\in \mathcal{M}$ and $h\in \{0,m\ell\}$ and let $a\in \mathcal{A}_{\ell,q}$ with $q>\max(H_i,H_j)+1/4$. There exists $\tau^2<+\infty$ such that the following convergence in distribution holds as $n\to +\infty$
\begin{equation}\label{eq-normCijm}
\sqrt{n} \big( C_{ij}^m(h) - \gamma_{ij}(h) \big) \stackrel{d}{\longrightarrow} \mathcal{N}(0,\tau^2).
\end{equation}
\end{lem}
We note that a central limit theorem also holds for $C_{ij}^m(-m\ell)$ since we recall that $C_{ij}^m(-m\ell)=C_{ji}^m(m\ell)$.
\begin{proof}
From the definition~(\ref{def-Cijmh}),
\begin{eqnarray*}
C_{ij}^m(h) - \gamma_{ij}^m(h) &=&  \frac1{n-m\ell-h} \sum_{k=m\ell+1}^{n-h} \left[ x_i^m(k)x_j^m(k+h)-\gamma_{ij}^m(h) \right] \\
&=& \frac1{n-m\ell-h} \sum_{k=1}^{n-m\ell-h} f(y(k)),
\end{eqnarray*}
where $y(k)=\left(x_i^m(k+m\ell),x_j(k+m\ell+h)\right)^t$ for $k=1,\ldots,n-m\ell-h$ and where 
$$
\begin{array}{rlcl}
f:&\R^2  &\to & \R \\
&y=(y_1,y_2)^t & \mapsto & y_1y_2 - \gamma_{ij}^m(h).
\end{array}
$$
Since $C_{ij}^m(h)-\gamma_{ij}^m(h)$ can be expressed as a centered empirical mean, the result is based on the application of a multivariate central limit theorem for non-linear functional of stationary Gaussian sequences obtained by \cite{Arco94}, Theorem~4. For this, note first that the Hermite rank  of the function $f$ (in the sense of \cite{Arco94}, Equation (2.2)) is two. Secondly, the condition on $q$ and Lemma~\ref{prop-asympgm1m2} (ii) (with $\alpha=2$) ensure that for any $i^\prime,j^\prime=i,j$
$$
\sum_{k\in \Z} \gamma_{i^\prime,j^\prime}^m(k+h)^2 =\sum_{k\in \Z} \gamma_{i^\prime,j^\prime}^m(k)^2 <+\infty.
$$
Theorem~4 of \cite{Arco94} can be applied, leading, as $n\to+\infty$, to the convergence in distribution of $\sqrt{n-m\ell-h}(C_{ij}^m(h) - \gamma_{ij}^m(h) )$ to a centered Gaussian random variable with finite variance (that we de not want  to explicit here). The result is obtained since $m$ and $h$ are fixed.
\end{proof}

\bigskip

Now, let us focus on the proof of Proposition~\ref{prop-norm}.

$(i)$ To prove this convergence, we follow the Cram\`er-Wold device \cite{Durr10} and prove that for any $\alpha \in \R^{D_{\mathcal{M},p}}$, $\alpha^t \sqrt{n} (C_n-\gamma)$ converges in distribution to $\alpha^t Z$ where $Z$ is a random normal vector. Let $\alpha\in \R^{D_{\mathcal{M},p}}$ be decomposed as follows
$$
\alpha=\left(
\alpha_{ii}^m,\alpha_{ij}^m,\alpha_{ij}^{m,+},\alpha_{ij}^{m,-}, \; (m\in \mathcal M,i,j=1,\ldots,p, j>i)
\right)^t,
$$
ordered as $C_n$ and $\gamma$. Then,
\begin{eqnarray*}
s_n &:=& \alpha^t (C_n -\gamma) \\
&=& \sum_i \sum_m \alpha_{ii}^m \; \frac1{n-m\ell} \sum_{k=m\ell+1}^n [x_i^m(k)^2-\gamma_{ii}^m(0)] \\
&&+\sum_{j>i} \sum_m \left\{ \alpha_{ij}^m \; \frac{1}{n-m\ell}\sum_{k=m\ell+1}^n [x_i^m(k)x_j^m(k)-\gamma_{ij}^m(0)] \right.\\
&&+ \alpha_{ij}^{m,+} \; \frac{1}{n-2m\ell}\sum_{k=m\ell+1}^{n-m\ell} [x_i^m(k)x_j^m(k+m\ell)-\gamma_{ij}^m(m\ell)]  \\
&&\left. +\alpha_{ij}^{m,-} \; \frac{1}{n-2m\ell}\sum_{k=m\ell+1}^{n-m\ell} [x_j^m(k)x_i^m(k+m\ell)-\gamma_{ij}^m(-m\ell)]  \right\}.
\end{eqnarray*}
Let $M = \max_{m\in \mathcal M}m$ ($M<+\infty$) and define
\begin{eqnarray*}
\widetilde{s}_n &:=& \frac1{n-2M\ell} \sum_{k=M\ell+1}^{n-M\ell} \bigg\{
  \sum_i \sum_m \alpha_{ii}^m [ x_i^m(k)^2 - \gamma_{ii}^m(0)] 
 \\
&&+\sum_{j>i}  \sum_m \bigg\{ \alpha_{ij}^m[x_i^m(k)x_j^m(k)-\gamma_{ij}(0)] + \alpha_{ij}^{m,+}[x_i^m(k)x_j^m(k+m\ell)-\gamma_{ij}(m\ell)]  \\
&& +\alpha_{ij}^{m,-}[x_j^m(k)x_i^m(k+m\ell)-\gamma_{ij}(-m\ell)] \bigg\} \bigg\}.
\end{eqnarray*}
Our first aim is to prove that $\sqrt{n}(s_n-\widetilde{s}_n) \stackrel{P}{\to} 0$ as $n\to +\infty$ (here and in the following $\stackrel{P}{\to}$ stands for the convergence in probability). For this, let us decompose the difference $s_n-\widetilde{s}_n=d_{1,n}+d_{2,n}$ where
\begin{eqnarray*}
d_{1,n} &:=& \sum_i \sum_m \frac1{n-2M\ell} \sum_{k\in I_m} \alpha_{ii}^m[x_i^m(k)^2-\gamma_{ii}^m(0)] \\
&&+\sum_{j>i} \sum_m \frac1{n-2M\ell} \sum_{k\in I_m} \alpha_{ij}^m[x_i^m(k)x_j^m(k)-\gamma_{ij}^m(0)] \\
&&+\sum_{j>i} \sum_{m\neq M} \frac1{n-2M\ell} \sum_{k\in I^{\pm}_m} \alpha_{ij}^{m,+}[x_i^m(k)x_j^m(k+m\ell)-\gamma_{ij}^m(m\ell)] \\
&&+\sum_{j>i} \sum_{m\neq M} \frac1{n-2M\ell} \sum_{k\in I^{\pm}_m} \alpha_{ij}^{m,-}[x_j^m(k)x_i^m(k+m\ell)-\gamma_{ij}^m(-m\ell)], 
\end{eqnarray*}
where for $m\neq M$, $I_m=\{m\ell+1,\ldots,M\ell\}\cup \{ n-M\ell,\ldots,n\}$,$I_m^\pm= \{m\ell+1,\ldots,M\ell\}\cup \{ n-M\ell,\ldots,n-m\ell\}$ and $I_M=\{n-M\ell,\ldots,n\}$. The remainder term $d_{2,n}$ is given by
\begin{eqnarray*}
d_{2,n} &:=&\frac{(m-2M)\ell}{n-2M\ell}  \left\{  \sum_i \sum_m  \alpha_{ii}^m C_{ii}^m(0)+\sum_{j>i} \sum_m \alpha_{ij}^m C_{ij}^m(0)\right\} \\
&&+\frac{2(m-M)\ell}{n-2M\ell}\sum_{j>i} \sum_{m\neq M} [\alpha_{ij}^{m,+} C_{ij}^m(m\ell) + \alpha_{ij}^{m,-} C_{ij}^m(-m\ell)].
\end{eqnarray*}
Let $e_{1,n}$ and $e_{2,n}$ be the firt (generic) sum terms of $d_{1,n}$ and $d_{2,n}$ given by
\begin{eqnarray*}
e_{1,n} &=& \frac1{n-2M\ell} \sum_{k\in I_m} \alpha_{ii}^m[x_i^m(k)^2-\gamma_{ii}^m(0)]\\
e_{2,n} &=& \frac{(m-2M)\ell}{n-2M\ell}  \alpha_{ii}^m C_{ii}^m(0).
\end{eqnarray*}
We now prove that $\sqrt{n}e_{k,n}\stackrel{P}{\to}0$ ($k=1,2$). The other terms follow similar arguments. Since the sum in $e_{1,n}$ contains a finite number of Gaussian random variables, $Var[\sqrt{n} e_{1,n}] = \mathcal{O}(n^{-1/2})$, which implies the convergence of $\sqrt{n}e_{1,n}$ to 0 in $L^2$ and so in probability. Finally, since $\sqrt{n}e_{2,n}= \alpha_{ii}\frac{(m-2M)\ell}{n-2M\ell} \sqrt{n}C_{ii}^m(0)$, Lemma~\ref{lem-normCijm} and Slutsky's Theorem ensure the expected convergence. 

As a consequence of the previous computations, we can concentrate ourselves on the asymptotic normality of $\widetilde s_n$. For this, let us define $y(k)=(x_i^m(k+m\ell+1),x_i^m(k+2m\ell+1), \; (m\in \mathcal{M}, i=1,\ldots,p))^t$ for $k=1,\ldots,n-2M\ell$ and
$$
\begin{array}{rlcl}
f_\alpha:&\R^{2|\mathcal{M}|p}& \to & \R \\
& y=(y_i^m,\widetilde{y}_i^m)^t & \mapsto& \sum_i \sum_m \alpha_{ii}^m[(y_i^m)^2-\gamma_{ii}^m(0)] + \sum_{j>i} \sum_m \alpha_{ij}^m[y_i^m y_j^m-\gamma_{ij}^m(0)] \\
&&& +\sum_{j>i} \sum_m \big( \alpha_{ij}^{m,+}[y_i^m \widetilde{y}_j^m-\gamma_{ij}^m(m\ell)] +   \alpha_{ij}^{m,-}[y_j^m \widetilde{y}_i^m-\gamma_{ij}^m(-m\ell)]\big).
\end{array}
$$
Then $\widetilde{s}_n$ can be expressed as the following empirical mean
$$
\widetilde{s}_n = \frac1{n-2M\ell} \sum_{k=1}^{n-2M\ell} f_\alpha(y(k)).
$$
For any vector $\alpha$, the Hermite rank of the function $f_\alpha$ is 2. Now, similarly as the proof of  Lemma~\ref{lem-normCijm}, we notice that for all $i,j=1,\ldots,p$, $m_1,m_2\in \mathcal{M}$, $h<+\infty$,
$$
\sum_{k\in \Z} \gamma_{ij}^{m_1,m_2}(k+h)^2 = \sum_{k\in \Z} \gamma_{ij}^{m_1,m_2}(k)^2 <+\infty,
$$
as soon as $q>H^\vee+1/4$. We can therefore apply Theorem 4 of \cite{Arco94} and obtain as $n\to +\infty$
$$
\sqrt{n} s_n \stackrel{d}{\longrightarrow} \mathcal{N} \left(0, \tau^2:=\sum_{k\in \Z} \E[ f_\alpha(y(t)) f_\alpha(y(t+k)) ]\right).
$$
The variance $\tau^2$ is given by $\tau^2 = \alpha^t \Sigma \alpha$. According to the definition of $\widetilde s_n$, the matrix $\Sigma$ can be partitioned into
\begin{eqnarray}
\Sigma  =\left(\begin{array}{cccc}
\Sigma_1 & \Sigma_{12}  & \Sigma_{13} & \Sigma_{14} \\
  \Sigma_{12}^t &  \Sigma_{2} &  \Sigma_{23} & \Sigma_{24} \\
 \Sigma_{13}^t & \Sigma_{23}^t &  \Sigma_{3} & \Sigma_{34} \\
 \Sigma_{14}^t & \Sigma_{24}^t  & \Sigma_{34}^t & \Sigma_{4}
 \end{array}\right). \label{eq-Sigma}
\end{eqnarray}
The matrix $\Sigma_1$ is the $|\mathcal{M}|p\times |\mathcal{M}|p$ covariance matrix of the vector containing the $x_i^m(t)^2$, $\Sigma_k$ (for $k=2,3,4$) are the $d_{\mathcal{M},p}\times d_{\mathcal{M},p}$ ($d_{\mathcal{M},p}=|\mathcal{M}|p(p-1)/2$) covariance matrices of the vectors containing the $x_i^m(t)x_j^m(t)$, the $x_i^m(t)x_j^m(t+m\ell)$ and the $x_j^m(t)x_i^m(t+m\ell)$ respectively. Other matrices are cross-covariances matrices (with dimension $d_{\mathcal{M},p}\times d_{\mathcal{M},p}$). Thus, the dimension of $\Sigma$ is $D_{\mathcal{M},p}\times D_{\mathcal{M},p}$ where $D_{\mathcal{M},p}:=|\mathcal{M}|p+3d_{\mathcal{M},p}= |\mathcal{M}|p(3p-1)/2$. Generic elements of these matrices can be evaluated using for example Isserlis Formula, \cite{Iss18}. The notation used hereafter follow the ordering of the vector $C_n$ (\ref{def-Cn}): the indices $i,j,i_1,j_1,i_2,j_2$ vary from 1 to $p$ (such that $j>i, j_1>i_1, j_2>i_2$), the indices $m_1,m_2$ vary in $\mathcal{M}$
\begin{eqnarray}
\begin{array}{ccl}
\left( \Sigma_1\right)_{i,j}^{m_1,m_2} &=& \sum_{k\in \Z}\E \left[
(x_i^{m_1}(t)^2-\gamma_{ii}^{m_1}(0))(x_j^{m_2}(t+k)^2-\gamma_{jj}^{m_2}(0))
\right] \\
&=& 2 \sum_{k\in \Z} \gamma_{ij}^{m_1,m_2}(k)^2, \\
\left( \Sigma_2\right)_{i_1,j_1,i_2,j_2}^{m_1,m_2} &=& \sum_{k\in \Z}\E \left[
(x_{i_1}^{m_1}(t)x_{j_1}^{m_1}(t)-\gamma_{i_1j_1}^{m_1}(0))(x_{i_2}^{m_2}(t+k)x_{j_2}^{m_2}(t+k)-\gamma_{i_2j_2}^{m_2}(0))
\right] \\
&=& \sum_{k\in \Z} [ \gamma_{i_1i_2}^{m_1,m_2}(k) \gamma_{j_1j_2}^{m_1,m_2}(k) + \gamma_{i_1j_2}^{m_1,m_2}(k)\gamma_{j_1i_2}^{m_1,m_2}(k)], \\
\left( \Sigma_3\right)_{i_1,j_1,i_2,j_2}^{m_1,m_2} &=& \sum_{k\in \Z}\E \left[
(x_{i_1}^{m_1}(t)x_{j_1}^{m_1}(t+m_1 \ell)-\gamma_{i_1j_1}^{m_1}(m_1\ell))(x_{i_2}^{m_2}(t+k)x_{j_2}^{m_2}(t+m_2\ell +k)-\gamma_{i_2j_2}^{m_2}(m_2\ell))
\right] \\
&=& \sum_{k\in \Z} [ \gamma_{i_1i_2}^{m_1,m_2}(k) \gamma_{j_1j_2}^{m_1,m_2}(k+(m_2-m_1)\ell) + \gamma_{i_1j_2}^{m_1,m_2}(k+m_2\ell)\gamma_{j_1i_2}^{m_1,m_2}(k-m_1\ell)], \\
\left( \Sigma_4\right)_{i_1,j_1,i_2,j_2}^{m_1,m_2} &=& \sum_{k\in \Z}\E \left[
(x_{j_1}^{m_1}(t)x_{i_1}^{m_1}(t+m_1 \ell)-\gamma_{i_1j_1}^{m_1}(-m_1\ell))(x_{j_2}^{m_2}(t+k)x_{i_2}^{m_2}(t+m_2\ell +k)-\gamma_{i_2j_2}^{m_2}(-m_2\ell))
\right] \\
&=&\left( \Sigma_3\right)_{j_1,i_1,j_2,i_2}^{m_1,m_2} \\
\left( \Sigma_{12} \right)_{i,i_1,j_1}^{m_1,m_2} &=& \sum_{k\in \Z} \E \left[ 
( x_i^{m_1}(t)^2-\gamma_{ii}^{m_1}(0)) ( x_{i_1}^{m_2}(t+k) x_{j_1}^{m_2}(t+k) -\gamma_{i_1j_1}^{m_2}(0))
\right] \\
&=& 2 \sum_{k\in \Z} \gamma_{ii_1}^{m_1,m_2}(k)\gamma_{ij_1}^{m_1,m_2}(k),\\
\left( \Sigma_{13} \right)_{i,i_1,j_1}^{m_1,m_2} &=& \sum_{k\in \Z} \E \left[ 
( x_i^{m_1}(t)^2-\gamma_{ii}^{m_1}(0)) ( x_{i_1}^{m_2}(t+k) x_{j_1}^{m_2}(t+k+m_2\ell) -\gamma_{i_1j_1}^{m_2}(m_2\ell))
\right] \\
&=& 2 \sum_{k\in \Z} \gamma_{ii_1}^{m_1,m_2}(k)\gamma_{ij_1}^{m_1,m_2}(k+m_2\ell),\\
\left( \Sigma_{14} \right)_{i,i_1,j_1}^{m_1,m_2} &=& \left( \Sigma_{13} \right)_{i,j_1,i_1}^{m_1,m_2}, \\
\left( \Sigma_{23}\right)_{i_1,j_1,i_2,j_2}^{m_1,m_2} &=& \sum_{k\in \Z}\E \left[
(x_{i_1}^{m_1}(t)x_{j_1}^{m_1}(t)-\gamma_{i_1j_1}^{m_1}(0))(x_{i_2}^{m_2}(t+k)x_{j_2}^{m_2}(t+m_2\ell +k)-\gamma_{i_2j_2}^{m_2}(m_2\ell))
\right] \\
&=& \sum_{k\in \Z} [ \gamma_{i_1i_2}^{m_1,m_2}(k) \gamma_{j_1j_2}^{m_1,m_2}(k+m_2\ell) + \gamma_{i_1j_2}^{m_1,m_2}(k+m_2\ell)\gamma_{j_1i_2}^{m_1,m_2}(k)], \\
\left( \Sigma_{24}\right)_{i_1,j_1,i_2,j_2}^{m_1,m_2} &=& \left( \Sigma_{23}\right)_{j_1,i_1,j_2,i_2}^{m_1,m_2}, \\ 
\left( \Sigma_{34}\right)_{i_1,j_1,i_2,j_2}^{m_1,m_2} &=& \sum_{k\in \Z}\E \left[
(x_{i_1}^{m_1}(t)x_{j_1}^{m_1}(t+m_1\ell)-\gamma_{i_1j_1}^{m_1}(m_1\ell))(x_{j_2}^{m_2}(t+k)x_{i_2}^{m_2}(t+m_2\ell +k)-\gamma_{i_2j_2}^{m_2}(-m_2\ell))
\right] \\
&=& \sum_{k\in \Z} [ \gamma_{i_1j_2}^{m_1,m_2}(k) \gamma_{j_1i_2}^{m_1,m_2}(k+(m_2-m_1)\ell) + \gamma_{i_1i_2}^{m_1,m_2}(k+m_2\ell)\gamma_{j_1j_2}^{m_1,m_2}(k-m_1\ell)] \\
&=& \left(\Sigma_3 \right)_{i_1,j_1,j_2,i_2}^{m_1,m_2}. \\
\end{array}
\label{matrice-cov-def:eq}
\end{eqnarray}

$(ii)$ Recall that the function $g:\R^{D_{\mathcal{M},p}}\to \R^{p(p+1)}$ maps the vector $C_n$ to $\widehat{\theta}$, i.e. $g(C_n)=\widehat{\theta}$. Now, we leave the reader to verify that replacing the vector $C_n$ by $\gamma$ allows us to retrieve $\theta$, i.e. $g(\gamma)=\theta$. Thus in view of applying the delta method \cite{LehmC98}, we only have to prove the differentiability of $g$ in $\gamma$. 
We do not want to provide all these heavy justifications and computations which are not very informative. We only focus on the terms that could lead to a problem that is the term related to $\rho_{ij}$ and $\eta_{ij}$. If $|\mathcal{M}|=1$, the estimates of $\rho_{ij}$ and $\eta_{ij}$ reduce to~(\ref{naturalRhoEta}) and it is simple to check that the function $g$ is differentiable in $\gamma$ for any $\gamma$, which means for any dilations set $\mathcal M$ and any set of parameters $H,\sigma^2,\rho,\eta$ (ensuring the model is well-defined) and in particular for some components of $\rho$ and/or $\eta$ set to zero. Due to the absolute values used in the definition of $\rho_{ij}$ and $\eta_{ij}$ (\ref{def-rhoijest},\ref{def-etaijest},\ref{def-rhoetaijest}) when $|\mathcal{M}|>1$, a problem of differentiability could appear. We show hereafter that it is not the case.
We just choose an example, namely the partial derivatives of $\widehat{\rho}_{ij}$ and assert that the other terms follow similar arguments. We make an abuse of notation in the following but we believe the context is clear. So, let us focus on the definition of $\widehat{\rho}_{ij}$ (\ref{def-rhoijest},\ref{def-rhoetaijest}) where again (only) for the sake of simplicity of the presentation we assume that $\widehat{\pi}^a_{ij}(0)={\pi}^a_{ij}(0)$. Then, let $j>i$
\begin{eqnarray*}
\frac{\partial \widehat{\rho}_{ij}}{\partial C_{ij}^m(0)} (C_n) &=& 
\frac{\partial |\widehat{\rho}_{ij} |}{\partial C_{ij}^m(0)} (C_n) \times \sign(C_{ij}^m(0))  \\
&=&\frac1{|\mathcal{M}|}\;  \sign(C_{ij}^m(0)) \frac{|\widehat{\rho}_{ij}|}{|C_{ij}^m(0)|}  \times \sign(C_{ij}^m(0))  \\
&=& \frac1{|\mathcal{M}|}\;  \frac{|\widehat{\rho}_{ij}|}{|C_{ij}^m(0)|}.
\end{eqnarray*}
And therefore, when evaluated at $\gamma$ we obtain
\begin{eqnarray*}
\frac{\partial \widehat{\rho}_{ij}}{\partial C_{ij}^m(0)} (\gamma) &=& \frac1{|\mathcal{M}|}\; \frac{|\rho_{ij}|}{m^{H_i+H_j} \sigma_i \sigma_j|\rho_{ij}| \pi^a_{ij}(0)} = \frac1{|\mathcal{M}|}\; \frac{1}{m^{H_i+H_j} \sigma_i \sigma_j \pi^a_{ij}(0)},
\end{eqnarray*}
which does not involve any continuity problem for the whole set of parameters $\mathcal{M},H,\sigma^2,\rho$ and $\eta$.
In the same spirit for example when differentiating $\widehat{\rho}_{ij}$ with respect to $C_{ii}^m(0)$ we get
\begin{eqnarray*}
\frac{\partial \widehat{\rho}_{ij}}{\partial C_{ii}^m(0)} (C_n) = -\frac1{2|\mathcal{M}|} \frac{|\widehat{\rho}_{ij}|}{C_{ii}^m(0)} \times  \sign(C_{ij}^m(0)),
\end{eqnarray*}
leading to 
$$
\frac{\partial \widehat{\rho}_{ij}}{\partial C_{ii}^m(0)} (\gamma) = -\frac{1}{2|\mathcal{M}|\gamma_{ii}^m(0)} |\rho_{ij}|\times \sign (\rho_{ij}) = -\frac{1}{2|\mathcal{M}|\gamma_{ii}^m(0)} \rho_{ij},
$$
and the conclusion is the same as previously.

\clearpage

\begin{center}
\begin{table}[ht]
\hspace*{0cm}
\begin{tabular}{|lr|ccc|ccc|ccc|}
\hline 
\multicolumn{2}{|c}{Parameters}& \multicolumn{3}{|c|}{causal mfBm} & \multicolumn{3}{c|}{well-balanced mfBm} & \multicolumn{3}{c|}{general mfBm ($\eta=0.2$)}\\  
  \cline{3-11}
 && $\widehat{H}^v$ & $\widehat{H}^c$ & $\widehat{H}^d$ & $\widehat{H}^v$ & $\widehat{H}^c$ & $\widehat{H}^d$ & $\widehat{H}^v$ & $\widehat{H}^c$ & $\widehat{H}^d$\\ 
  \hline
 \textbf{H=0.2}& $\rho=0.1$ & 0.0006 & 0.0078 & 0.0126 & 0.0005 & 0.0062 & 0.0203 & 0.0005 & 0.0060 & 0.0131 \\ 
 $(p=2)$ &$0.5$    & 0.0005 & 0.0007 & 0.0140 & 0.0005 & 0.0006 & 0.0119 & 0.0006 & 0.0008 & 0.0126 \\ 
 &$0.9$   & 0.0005 & 0.0006 & 0.0140 & 0.0005 & 0.0005 & 0.0148 & 0.0005 & 0.0006 & 0.0149 \\ 
 \hline
 \textbf{H=0.5}& $\rho=0.1$ & 0.0009 & 0.0136 & 0.0171 & 0.0007 & 0.0100 & 0.0204 & 0.0010 & 0.0085 & 0.0111 \\ 
  $(p=2)$ &$0.5$  & 0.0008 & 0.0011 & 0.0132 & 0.0008 & 0.0011 & 0.0169 & 0.0008 & 0.0012 & 0.0125 \\ 
 &$0.9$   & 0.0008 & 0.0009 & 0.0140 & 0.0009 & 0.0010 & 0.0146 & 0.0010 & 0.0010 & 0.0126 \\ 
  \hline
\textbf{H=0.8} &$\rho=0.1$   & 0.0010 & 0.0085 & 0.0148 & 0.0009 & 0.0110 & 0.0176 & 0.0009 & 0.0157 & 0.0194 \\ 
 $(p=2)$ &$0.5$  & 0.0008 & 0.0012 & 0.0099 & 0.0012 & 0.0015 & 0.0114 & 0.0010 & 0.0014 & 0.0114 \\ 
 &$0.9$  & 0.0009 & 0.0009 & 0.0104 & 0.0009 & 0.0009 & 0.0126 & 0.0009 & 0.0010 & 0.0126 \\ 
  \hline
 \textbf{H=0.1:0.5}& $\rho=0.1$& 0.0008 & 0.0048 & 0.0143 & 0.0007 & 0.0054 & 0.0155 & 0.0006 & 0.0052 & 0.0139 \\ 
  $(p=2)$ &$0.5$  & 0.0006 & 0.0007 & 0.0116 & 0.0007 & 0.0010 & 0.0113 & 0.0007 & 0.0010 & 0.0154 \\ 
 &$0.9$  & $\times$ & $\times$  & $\times$  & $\times$  & $\times$  & $\times$  & $\times$  & $\times$  & $\times$  \\
  \hline
 \textbf{H=0.5:0.9} $\quad$ & $\rho=0.1$ & 0.0009 & 0.0125 & 0.0116 & 0.0009 & 0.0073 & 0.0125 & 0.0008 & 0.0035 & 0.0135 \\ 
  $(p=2)$&$0.5$ & 0.0008 & 0.0009 & 0.0090 & 0.0009 & 0.0011 & 0.0099 & 0.0009 & 0.0011 & 0.0112 \\ 
 &$0.9$ & $\times$  & $\times$  & $\times$  & $\times$  & $\times$  & $\times$  &  $\times$ & $\times$  &  $\times$ \\ 
  \hline
 \textbf{H=0.2}& $\rho=0.1$& 0.0005 & 0.0130 & 0.0255 & 0.0005 & 0.0135 & 0.0270 & 0.0006 & 0.0141 & 0.0243 \\ 
 $(p=5)$ &$0.5$  & 0.0006 & 0.0011 & 0.0234 & 0.0005 & 0.0011 & 0.0225 & 0.0006 & 0.0010 & 0.0226 \\ 
&$0.9$    & 0.0006 & 0.0006 & 0.0232 & 0.0005 & 0.0006 & 0.0250 & 0.0006 & 0.0007 & 0.0248 \\ 
  \hline
 \textbf{H=0.8}& $\rho=0.1$& 0.0010 & 0.03050 & 0.0265 & 0.0008 & 0.0292 & 0.0239 & 0.0009 & 0.0255 & 0.0249 \\ 
 $(p=5)$& $0.5$  & 0.0009 & 0.0017 & 0.0202 & 0.0011 & 0.0021 & 0.0185 & 0.0011 & 0.0021 & 0.0198 \\ 
&$0.9$  & 0.0009& 0.0010& 0.0199 & 0.0011& 0.0013 & 0.0214 & 0.0010 & 0.0011 & 0.0192 \\ 
  \hline
 \textbf{H=0.1:0.5}& $\rho=0.1$& 0.0006 & 0.0163 & 0.0272 & 0.0006 & 0.0172 & 0.0283 & 0.0006 & 0.0143 & 0.0251 \\ 
  $(p=5)$ &$0.5$ & 0.0007 & 0.0014 & 0.0220 & 0.0006 & 0.0012 & 0.0247 & 0.0006 & 0.0013 & 0.0229 \\ 
 &$0.9$  & $\times$  & $\times$  & $\times$  &  $\times$ & $\times$  & $\times$  & $\times$  &  $\times$ & $\times$  \\ 
 \hline
   \end{tabular}
\caption{Empirical means of MSE of $H$ estimates based on 100 replications of causal, well-balanced and general mfBm of length $n=1000$. 
For the general case, the parameter $\eta_{ij}$ is fixed to $.2\times (1-H_i-H_j)$ for $i>j$ and $\eta_{ji}=-\eta_{ij}$. The letters $v,c,d$ 
respectively correspond to the Hurst exponents estimators computed with the weights $w=(1,0,0)$ (using only the empirical variances) , $w=(1,1,0)$ 
(using in addition the empirical covariances) and $w=(1,1,1)$ (using in addition the difference of empirical covariances at lags $\pm m$).}
\label{hest:tb}
\end{table}
\end{center}
  \begin{center}
\begin{table}[ht]
\begin{tabular}{|lr|cc|cc|cc|}
\hline 
\multicolumn{2}{|c}{Parameters}& \multicolumn{2}{|c|}{causal mfBm} & \multicolumn{2}{c|}{well-balanced mfBm} & \multicolumn{2}{c|}{general mfBm ($\eta=0.2$)}\\  
  \cline{3-8}
 & & $\widehat{\rho}^v$ & $\widehat{\rho}^d$ &
   $\widehat{\rho}^v$ & $\widehat{\rho}^d$ &
   $\widehat{\rho}^v$ & $\widehat{\rho}^d$  \\ 
  \hline
 \textbf{H=0.2}& $\rho=0.1$ &  0.0013 & 0.0013 &   0.0012 &  0.0078 & 0.0014 &0.0078 \\
  $(p=2)$ &$0.5$    &  0.0007 & 0.0007  & 0.0007 &  0.1878 & 0.0007 &  0.1878\\
 &$0.9$   &  0.0001 & 0.0001  & 0.0001 &  0.6080 & 0.0001&  0.6080 \\
 \hline
\textbf{H=0.8} &$\rho=0.1$   &  0.0017 & 0.0017  & 0.0022&  0.0080  & 0.0022 & 0.0081\\ 
 $(p=2)$ &$0.5$   & 0.0011&  0.0041 & 0.0011 &  0.1883  & 0.0013&  0.1890 \\ 
 &$0.9$   & 0.0001&   0.0007  & 0.0001&  0.6084 & 0.0001&  0.6083  \\
 \hline
 \textbf{H=0.1:0.5}& $\rho=0.1$ & 0.0014 &  0.0013 & 0.0010&  0.0078  & 0.0009 &  0.0077\\ 
  $(p=2)$ &$0.5$   & 0.0005 &  0.0007 & 0.0005&  0.1878  & 0.0007 &  0.1879\\ 
 &$0.9$  & $\times$  & $\times$  & $\times$  & $\times$  & $\times$  & $\times$  \\
  \hline
 \textbf{H=0.5:0.9} $\quad$ & $\rho=0.1$  & 0.0012 & 0.0027 & 0.0010 &0.0082 & 0.0011&  0.0083 \\
  $(p=2)$&$0.5$  & 0.0037&   0.0411 & 0.0023 &  0.2050 & 0.0028 &  0.2029 \\
 &$0.9$ & $\times$   & $\times$  & $\times$  &  $\times$ & $\times$  &  $\times$ \\ 
  \hline
 \textbf{H=0.2}& $\rho=0.1$ & 0.0011&  0.0011 & 0.0012 &  0.0065  & 0.0013&  0.0065 \\
 $(p=5)$ &$0.5$   & 0.0007&  0.0008 & 0.0006 &  0.1503  & 0.0008 &  0.1504\\ 
&$0.9$     & 0.0001 &   0.0003& 0.0001 &  0.4869  & 0.0001&  0.4864 \\ 
  \hline
 \textbf{H=0.8}& $\rho=0.1$ & 0.0019&  0.0039  & 0.0019&  0.0075  & 0.0018 &  0.0076\\ 
 $(p=5)$& $0.5$   & 0.0010 &  0.0686 & 0.0010 &  0.1753 & 0.0012 &  0.1765\\
&$0.9$   & 0.0001&  0.2060  & 0.0001 &  0.5683 & 0.0001 &  0.5693\\ 
  \hline
 \textbf{H=0.1:0.5}& $\rho=0.1$& 0.0013 & 0.0013  & 0.0013 &  0.0065& 0.0013 &  0.0065 \\
  $(p=5)$ &$0.5$  & 0.0006 &  0.0010 & 0.0007 &  0.1514 & 0.0007 &  0.1506\\ 
 &$0.9$  & $\times$  &  $\times$  & $\times$  & $\times$  &  $\times$ & $\times$  \\ 
 \hline
   \end{tabular}
\caption{Empirical means of MSE of the estimates of $\rho_{ij}$ based on 100 replications of causal, well-balanced and general mfBm of length $n=1000$. 
For the general case, the parameter $\eta_{ij}$ is fixed to $.2\times (1-H_i-H_j)$ for $i>j$ and $\eta_{ji}=-\eta_{ij}$. The letters $n,d,v$ correspond to different strategies detailed in Section \ref{experiment:sec}}
\label{rhoest:tb}
\end{table}
   \end{center}
  \begin{center}
\begin{table}[ht]
\begin{tabular}{|lr|cc|cc|cc|}
\hline 
\multicolumn{2}{|c}{Parameters}& \multicolumn{2}{|c|}{causal mfBm} & \multicolumn{2}{c|}{well-balanced mfBm} & \multicolumn{2}{c|}{general mfBm ($\eta=0.2$)}\\  
  \cline{3-8}
 && $\widehat{\eta}^v$ & $\widehat{\eta}^d$ &
   $\widehat{\eta}^v$ & $\widehat{\eta}^d$ &
   $\widehat{\eta}^v$ & $\widehat{\eta}^d$  \\ 
  \hline
 \textbf{H=0.2}& $\rho=0.1$ &   0.2081 &0.2663   & 0.1834&  0.3598  & 0.1169 &  3.8969\\
 $(p=2)$ &$0.5$     & 0.1513&  0.4219  & 0.1315 & 163.75 & 0.0862 &  0.3576\\ 
 &$0.9$    & 0.0401 &  1.3426 & 0.0366 &  0.0521 & 0.0112 &  0.7055\\
 \hline
\textbf{H=0.8} &$\rho=0.1$   & 0.1698  &  15.6143 & 0.1240 & 0.1985  & 0.0908 & 0.4624\\
 $(p=2)$ &$0.5$   & 0.1181&  0.1155  & 0.1058 & 0.1124 & 0.0545 &  0.2181\\
 &$0.9$   & 0.0274 &   0.0312 & 0.0353 & 0.0609 & 0.0065 & 0.0797 \\
  \hline
 \textbf{H=0.1:0.5}& $\rho=0.1$ & 0.2416 & 26.3081  & 0.3332 & 4.6943  & 0.2333 & 8.3665\\
  $(p=2)$ &$0.5$   & 0.0507 &    0.4746  & 0.2185  & 684.53& 0.1839 & 2.0666\\
 &$0.9$  & $\times$ &  $\times$  & $\times$  & $\times$  & $\times$  & $\times$  \\
  \hline
 \textbf{H=0.5:0.9} $\quad$ & $\rho=0.1$  & 0.2472 &   1.6664  & 0.2067 & 1.6916 & 0.0999 & 1.6900\\
  $(p=2)$&$0.5$   & 0.9107 & 15.6106  & 0.0927 & 0.9720 & 0.0519  & 0.6583\\ 
 &$0.9$   & $\times$  & $\times$  & $\times$  &  $\times$ & $\times$  &  $\times$ \\ 
  \hline
 \textbf{H=0.2}& $\rho=0.1$  & 0.1919  & 0.8977 & 0.1954  & 16.3425 & 0.1095 & 0.4842\\
 $(p=5)$ &$0.5$   & 0.1440 &     0.9434 & 0.1501  & 0.6351 & 0.0790 & 32.3723 \\ 
&$0.9$      & 0.0362 & 1.3560  & 0.0394 & 0.2708  & 0.0096 & 0.4003\\ 
  \hline
 \textbf{H=0.8}& $\rho=0.1$  & 0.1418 & 3.8831  & 0.1503 & 0.3738  & 0.0789 & 0.8545\\
 $(p=5)$& $0.5$   & 0.1126& 0.3640   & 0.1154 & 0.2667  & 0.0522 & 0.2319\\
&$0.9$    & 0.0306 & 0.1154  & 0.0304 & 0.0547  & 0.0074 & 0.0494\\ 
  \hline
 \textbf{H=0.1:0.5}& $\rho=0.1$ & 2.8168&900.01   &  1.6256 & 96.4169 & 2.5045 & 712.66 \\
  $(p=5)$ &$0.5$  & 30.279 &  259.83&  & 2.7290 $\times$ & 14.264 &  1455.58\\ 
 &$0.9$  &   $\times$ & $\times$  & $\times$  & $\times$  &  $\times$ & $\times$  \\ 
 \hline
   \end{tabular}
\caption{Empirical means of MSE of the estimates of $\eta_{ij}$ based on 100 replications of causal, well-balanced and general mfBm of length $n=1000$. 
For the general case, the parameter $\eta_{ij}$ is fixed to $.2\times (1-H_i-H_j)$ for $i>j$ and $\eta_{ji}=-\eta_{ij}$. The letters $n,d,v$ correspond to different strategies detailed in Section \ref{experiment:sec}}
\label{etaest:tb}
\end{table}
\end{center}


\bibliographystyle{apalike}

\end{document}